
\input epsf
\input amstex
\magnification=1200
\documentstyle{amsppt}
\NoBlackBoxes
\def\1{\bold 1}
\def\a{\alpha}
\def\l{\lambda}
\def\lk{{\lambda^k}}
\def\g{\frak g}
\def\h{\frak h}
\def\Z{\Bbb Z}
\def\C{\Bbb C}
\def\R{\Bbb R}
\def\K{\Bbb K}

\def\i{\text{\rm i}}
\def\Tr{\operatorname{Tr}}
\def\Id{\operatorname{Id}}
\def\ch{\operatorname{ch }}
\def\dim{\operatorname{dim }}
\def\Hom{\operatorname{Hom}}
\def\End{\operatorname{End}}
\def\What{{\widetilde W}}
\def\v{^\vee}
\def\<{\langle}
\def\>{\rangle}

\def\sltwo{\frak s\frak l _2 }
\def\sln{\frak{sl}_n}
\def\Ug{U_q\frak g}

\def\o{\otimes}
\def\eps{\varepsilon}

\topmatter

\title On inner product in modular tensor categories. I
\endtitle
\author
Alexander A. Kirillov, Jr.\\
\it Dept. of Mathematics, MIT\\
Cambridge, MA 02139, USA\\
 e-mail: kirillov\@math.mit.edu\\
http://web.mit.edu/kirillov/www/home.html
\endauthor
\leftheadtext{Alexander Kirillov, Jr.}



\date November 15, 1995\enddate

\endtopmatter

\document

\centerline{\bf q-alg/9508017}

\head 0. Introduction\endhead

In this paper we study some properties of tensor categories that arise in
2-dimensional conformal and 3-dimensional topological quantum field
theory -- so called modular tensor categories. By definition, these
categories are braided tensor categories with duality which are
semisimple, have finite number of simple objects and satisfy some
non-degeneracy condition. Our main example of
such a category is the reduced category of representations of a
quantum group $\Ug$ in the case when $q$ is a root of unity (see
\cite{AP, GK}).

The main property of such categories is that we can introduce a
natural projective action of mapping class group of any 2-dimensional
surface with marked points on appropriate spaces of morphisms in this
category (see \cite{Tu}). This property explains the name ``modular
tensor category'' and is crucial for establishing relation with
3-dimensional quantum field theory and in particular, for construction
of invariants of 3-manifolds (Reshetikhin-Turaev invariants).

 In particular, for the torus with one puncture we get a projective
action of the modular group $SL_2(\Z)$ on any space of morphisms $\Hom
(H, U)$, where $U$ is any simple object and $H$ is a special object
which is an analogue of regular representation (see \cite{Lyu}).  In
the case $U=\C$ this action is well known: it is the action of modular
group on the characters of corresponding affine Lie algebra. We study
this action for arbitrary representation $U$; in particular, we show
that this action is unitary with respect to a natural inner product on
the space of intertwining operators.

In the special case $\g=\sln$ and $U$ being a symmetric power of
fundamental representation this is closely related with Macdonald's
theory. It was shown in the paper
\cite{EK3} (though we didn't use the word ``$S$-matrix'' there)
 that in this case the matrix coefficients of the matrix $S$
are some special values of Macdonald's polynomials of type $A_{n-1}$.
Thus,  the properties of $S$-matrix immediately yield a
number of identities for values of Macdonald's polynomials at
roots of $1$. In this case, the action of modular group is closely
related with the difference Fourier transform defined in a recent paper of
Cherednik (\cite{Ch}). In particular, this shows that for $\g=\sltwo$
all matrix elements of $S$-matrix can be written in terms of
$q$-ultraspherical polynomials.

Unfortunately, we had to spend a large part of this paper   recalling known
facts about modular tensor categories and quantum groups at roots of
unity; though these results are well-known to experts, they are
scattered in numerous papers, and some parts are not written anywhere
at all. Thus, Sections~1 and 3 and large part of Section~6 are
expository.

The paper is organized as follows. In Section~1 we recall basic facts
about modular tensor categories (MTC), in particular, the action of modular
group and various symmetries of this action. In Section~2 we define an
inner product on the space of intertwiners in modular tensor
categories with some additional properties (hermitian MTC's), and
prove that the action of modular group is unitary with respect to this
inner product.

In Section~3 we recall, following Andersen, construction of MTC from
representations of quantum groups at roots of unity. In Section~4 we
show that this category can be endowed with a natural hermitian
structure.

Section~5 is devoted to a special case of the constructions above;
namely, we let $\g=\sln$ and take $U$ to be a symmetric power of
fundamental representation. We show that in this case $S$-matrix can be
written in terms of values of Macdonald's polynomials of type $A_{n-1}$ at
roots of unity, which gives many identities for these special values.
These expressions coincide with Cherednik's formulas for difference
Fourier transform.

Sections~6 and 7 are devoted to further study of MTC's coming from
quantum groups at roots of unity. In particular, we describe the
Grothendieck ring of these categories (which is not new); we also
give another description of the hermitian structure on them.

In the next papers we will  apply the same construction to
the modular tensor category arising from the affine Lie algebras, in which
case it will give a modular invariant inner product on the space of
conformal blocks.

\subhead Acknowledgments \endsubhead

I'd like to  expresses my deep gratitude to my advisor Igor Frenkel
for his guidance and encouragement. In particular, it is from his
course that I first learned the structure of modular tensor category
related to representations of quantum groups, as well as many other
structures appearing in this paper.

Special thanks are due to Pavel Etingof. This paper continues the
ideas introduced in a series of our joint papers, and importance of
discussions with Pavel for this paper can not be overestimated.

Also, I'd like to thank Ivan Cherednik, Thomas Kerler, David Kazhdan,
George  Lusztig and Stephen Sawin for stimulating discussions and Harvard
Mathematics department for its hospitality during my work.

Financial support was provided by Alfred P. Sloan dissertation
fellowship.

\head 1. Modular tensor categories.\endhead

In this section we review the main definitions relating to modular
tensor categories. This section is completely expository and does not
contain any new results. We start with a quick introduction to the
notion of a ribbon category, introduced by Reshetikhin and Turaev
(\cite{RT1, RT2}); we refer the reader to recent books by Kassel
(\cite{Kas}) and Turaev (\cite{Tu}) for detailed exposition.

\subhead Ribbon categories and graphs\endsubhead

A {\it ribbon category} is an additive category $\Cal C$ with the
following additional structures:

\roster
\item
 A bifunctor $\o:\Cal C\times \Cal C\to \Cal C$ along with
 functorial associativity and commutativity isomorphisms:

$$\gathered
a_{V_1,V_2,V_3}:(V_1\o V_2)\o V_3\to V_1\o(V_2\o V_3),\\
\check R_{V,W}:V\o W\to W\o V.\endgathered$$

\item A unit object $\1\in \text{Obj }\Cal C$ along with isomorphisms
$\1\o V\to V, V\o \1\to V$.

\item A notion of dual: for every object $V$ we have a (left) dual
$V^*$ and homomorphisms

$$\gathered
e_V: V^*\o V\to \1, \\
i_V:\1\to V\o V^*\endgathered$$

\item Balancing, or a system of twists, i.e. functorial isomorphisms
$\theta_V:V\to V$, satisfying the compatibility condition
$$\theta_{V\o W}= \check R_{W,V} \check R_{V, W}(\theta_V\o
\theta_W).$$

\endroster

These structures have to obey a number of properties, the list of
which can be found in \cite{Kas}. Using them, one can define
functorial isomorphisms $\delta_V: V\to
V^{**}$ which is compatible with tensor products and unit object.
This, in particular, implies that for every $V$ we also have
the right dual ${}^*V$ which is canonically isomorphic to the left
dual and homomorphisms

$$\gathered
V\o {}^*V\to \1,\\
\1\to {}^*V\o V.\endgathered$$

In another terminology, ribbon categories are called braided monoidal
rigid balanced categories (these words refer to the data we introduced
in items (1)--(4) above, respectively).

Unless otherwise specified, we will assume that our category is abelian.
 We will also use the following theorem, due to MacLane:
each ribbon category is equivalent to a strict one, i.e. such a
category in which $(V_1\o V_2)\o V_3=V_1\o (V_2\o V_3)$ (not only
isomorphic but is the same object!), and associativity morphism is the
identity morphism; proof of this fact can be found in \cite{Mac}.
 Unless otherwise specified, we only consider strict categories, and
thus we can write tensor products of many objects without bothering
about the parentheses.

Ribbon tensor categories admit a nice pictorial representation: if
we have a directed tangle with braids labeled by objects of $\Cal C$
and coupons labeled by morphisms then we can assign to such a tangle a
morphism in category $\Cal  C$ by certain rules -- see \cite{RT1, RT2, Tu} or
\cite{Kas}. The theorem proved by
Reshetikhin and Turaev says that this morphism only depends on the
isotopy class of the tangle; thus, we can prove identities about
morphisms by manipulating with corresponding tangles.
  Also, note that if we replace a label $V$ of a
certain braid by $V^*$ and reverse the direction of this braid then we
get the same morphism (up to canonical isomorphisms  $ V\simeq
V^{**}$).

For technical reasons, we will draw lines instead of ribbons; the only
problem with that is that when establishing isotopy of graphs one must
be careful to count the  twists. Examples of tangles and
corresponding operators and some identities are shown on Figure~1.
Note that the operators act ``from bottom to top'', even though the
arrows are oriented downwards.

\midinsert
\centerline{\epsfbox{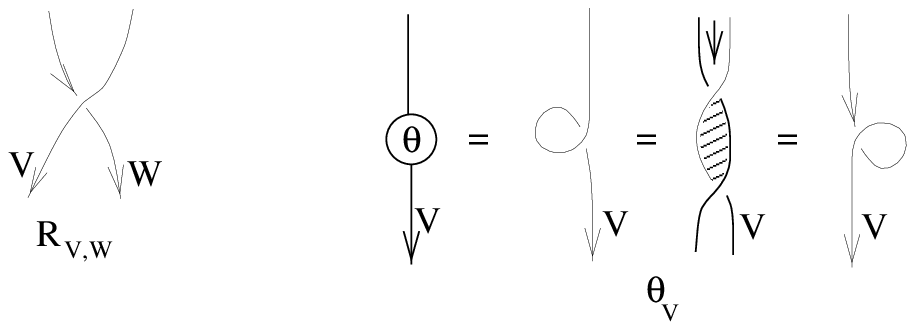}}
\botcaption{Figure 1a}
\endcaption
\endinsert

\midinsert
\centerline{\epsfbox{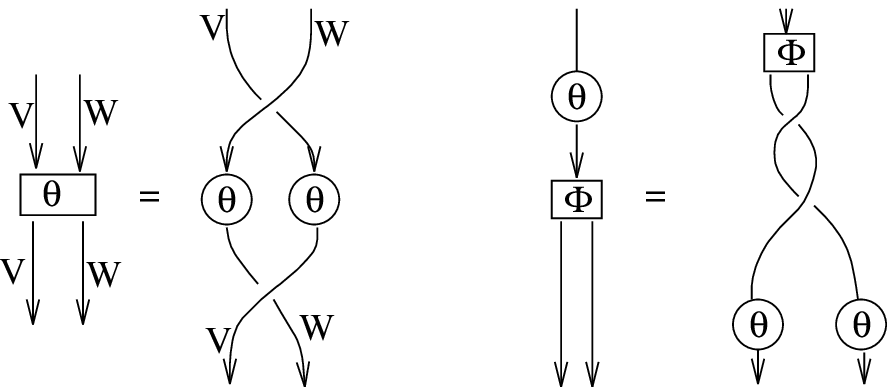}}
\botcaption{Figure 1b}
\endcaption
\endinsert

Let us additionally assume that our category is {\it semisimple},
i.e.

\roster
\item It is defined over some field $\K$: all the spaces of
homomorphisms are finite-dimensional vector spaces over $\K$.

\item Isomorphism classes of simple  objects in $\Cal C$
are indexed by elements of some
set $I$; we will use the notation $X_i, i\in I$ for the corresponding
simple, choosing the indexing so that $X_0=\1$. This
implies that we have an involution ${}^*:I\to I$ such that
$$X_i^*\simeq X_{i^*};$$
in particular, $0^*=0$.

\item ``Schur's Lemma'':

$$\Hom (X_i, X_j)=\cases \K, \quad i=j\\ 0,\quad i\ne j\endcases$$

\item Every object is completely reducible: every $V\in \text{Obj
}\Cal C$ can be written  in the form
$$V=\bigoplus_{i\in I} N_i X_i,$$
where $N_i\in \Z_+$, and the sum is finite (i.e., almost all
$N_i=0$).
\endroster

\remark{Remark} In fact, these axioms are abundant: for example, $\Hom
(X_i, X_j)=0$ for $i\ne j$ can be deduced from other axioms, see
\cite{Tu}.
\endremark

It will be convenient in the future to fix isomorphisms $X_{i^*}\simeq
(X_i)^*$ so that the composition

$$X_i=X_{i^{**}}\simeq (X_{i^*})^*\simeq X_i^{**}\tag 1.1$$
coincides with the map  $\delta_{X_i}$. This is equivalent to choosing
a nonzero homomorphism $X_{i^*}\o X_i\to \1$ (``Shapovalov form'').

Semisimplicity is a very restrictive requirement; it  implies a lot
of properties. For example, we can define the multiplicity
coefficients $N_{ij}^k$ by

$$X_i\o X_j=\bigoplus N_{ij}^k X_k,$$
then we have the following obvious properties:

$$\gathered
N_{ij}^k=\dim \Hom (X_i\o X_j, X_k)=\dim\Hom(X_i\o X_j\o X_k^*, \1),\\
N_{ij}^k=N_{ji}^k= N_{ik^*}^{j^*}= N_{i^* j^*}^{k^*},\\
N_{ij}^0=\delta_{i j^*}.\endgathered\tag 1.2$$

For  an object $V\in \text{Obj }\Cal C$ define its dimension $\dim
V\in \K$ by the following  picture:

$$\dim V=\quad\vcenter{\epsfbox{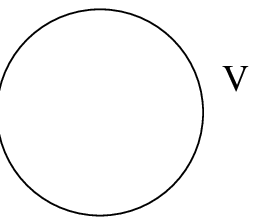}}\tag 1.3$$

More generally, for a morphism $f\in \Hom (V, V)$ we define its trace
$\Tr f$ by the following picture
$$\Tr f= \quad\vcenter{\epsfbox{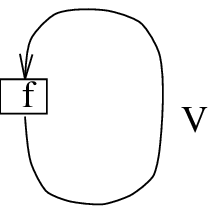}}$$

When objects of category are vector spaces over some field, the
dimension and trace defined above are usually called ``quantum
dimension'' (respectively, ``quantum trace'') to distinguish from
ordinary dimension and trace.

\proclaim{Lemma 1.1} \roster \item $\dim V^*=\dim V, \dim \1=1$.

\item $\dim V\o W=\dim V\cdot \dim W.$
\endroster\endproclaim

\subhead Action of modular group\endsubhead

As before, we assume that we have a semisimple ribbon category $\Cal
C$ with simple objects $X_i, i\in I$. Define the numbers $s_{ij}\in \K$ by the
following picture:

$$s_{ij}=\quad\vcenter{\epsfbox{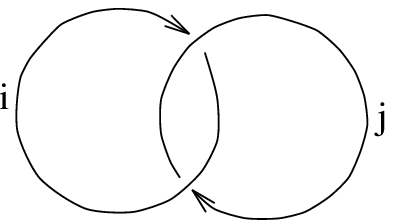}}\tag 1.4$$

(From now on, we will often label strands of tangles by the indices
$i\in I$ meaning by this $X_i$).

\proclaim{Proposition 1.2}

$$\gathered
s_{ij}=s_{ji}=s_{i^*j^*}=s_{j^*i^*},\\
s_{i0}=\dim X_i.\endgathered\tag 1.5$$
\endproclaim

\definition{Definition 1.3} A semisimple
ribbon category $\Cal C$ is called {\it modular} if it satisfies the
following properties:
\roster
\item It has only finite number of simple objects: $|I|<\infty$.

\item The matrix $s=(s_{ij})_{i,j\in I}$, where $s_{ij}$ is defined by
\rom{(1.4)},  is invertible.
\endroster
\enddefinition

We will give an example of a modular category later.

\remark{Remark} In fact, many authors (for example, Turaev) impose
weaker conditions, not necessarily requiring semisimplicity in our
sense. We are only
interested in the simplest case; thus the above definition is
absolutely sufficient for our purposes.  We refer the reader to
\cite{Ke} for discussion of non-semisimple case.
\endremark

\proclaim{Proposition 1.4} In a modular category, we have $\dim X_i\ne
0$ and

$$\epsfbox{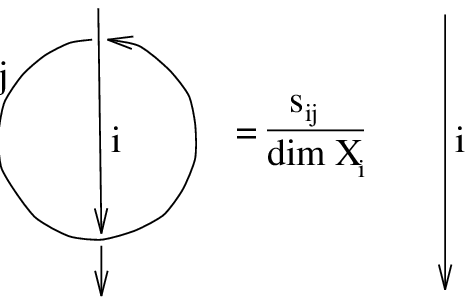}\tag 1.6$$

\endproclaim
The name ``modular'' is justified by the fact that in this case we can
define a projective action of the modular group $SL_2(\Z)$ on certain
objects in our category, which we will show below. To the best of my
knowledge, this construction
first appeared (in rather vague terms) in papers of Moore and Seiberg
(\cite{MS2,4}); later it was formalized by Lyubashenko (\cite{Lyu})
and others. Our exposition follows the book of Turaev.

The appearance of modular group in tensor categories may seem
mysterious; however, there is a simple geometrical explanation, based
on the fact that to each modular tensor category one can associate a
2+1-dimensional Topological Quantum Field Theory. This also shows that
in fact we have an action of mapping class group of any closed
oriented 2-dimensional surface on the appropriate objects in MTC. This
is the key idea of the book \cite{Tu}.

{}From now on, let us adopt the following convention: if some (closed)
strand on a picture is left unlabeled then we assume summation over
all labels $i\in I$ each taken with the weight $\dim X_i$. Then we
have the following propositions, proof of which (not too difficult)
can be found in \cite{Tu}.

Let us define
the numbers $\theta_i\in \K$ by

$$\theta_{X_i}=\theta_i \Id_{X_i},$$
then it is easy to see that $\theta_i=\theta_{i^*},
\theta_0=1$.

\proclaim{Proposition 1.5}

We have the following identities:

$$\epsfbox{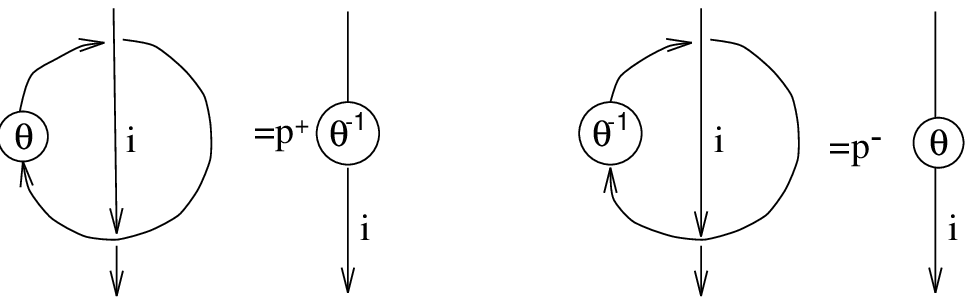},\tag 1.7$$
where
$$p^{\pm} = \sum_{i\in I} \theta_i^{\pm 1} (\dim X_i)^2.\tag 1.8$$

Also,

$$\vcenter{\epsfbox{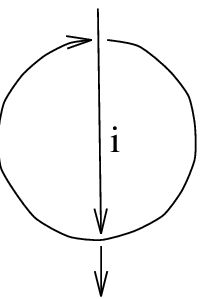}}\quad =p^+ p^-\delta_{i,0}.\tag 1.9$$
\endproclaim

\proclaim{Corollary 1.6}
\roster\item

$$\vcenter{\epsfbox{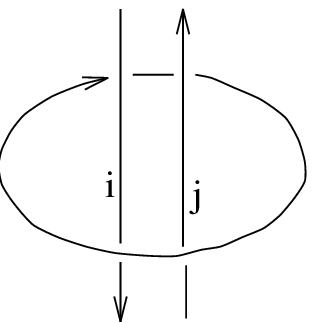}}\quad =p^+p^- \frac{\delta_{ij}}{\dim
X_i} \vcenter{\epsfbox{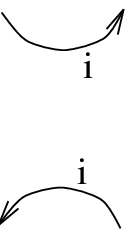}}
\tag 1.10$$
\item

$$p^{+}p^-=\sum (\dim X_i)^2.\tag 1.11$$

\endroster
\endproclaim

Define the matrices $t=(t_{ij})$ and $c=(c_{ij})$ (``charge conjugation
matrix'') by

$$\aligned
t_{ij}&=\delta_{ij}\theta_i,\\
c_{ij}&=\delta_{ij^*}.\endaligned
\tag 1.12$$
 Then it is easy to
deduce from Proposition~1.5 the following theorem:

\proclaim{Theorem 1.7} The matrices $s,t$ defined above satisfy the
following relations:

$$\gathered
 s^2=p^+p^- c,\\
( st)^3=p^+ s^2,\\
 s^2 t=t  s^2,\endgathered\tag 1.13$$
where $p^\pm$ are defined by \rom{(1.8)}.
\endproclaim

It is convenient to renormalize these matrices. Namely, let
 us assume that the following fractional powers exist in $\K$:

$$\gathered
D=\sqrt{p^+p^-}=\sqrt{\sum_{i\in I}(\dim X_i)^2},\\
\zeta= (p^+/p^-)^{1/6}\endgathered\tag 1.14$$
(we choose the roots so that $D\zeta^3=p^+$).
It follows from non-degeneracy of $s$ that $D,\zeta\ne 0$.

Define renormalized matrices

$$\tilde s=\frac{s}{D},\qquad \tilde t= \frac{t}{\zeta}.\tag 1.15$$
Then Theorem 1.7 is rewritten in the following form:

$$\gathered
\tilde s^2=c,\\
(\tilde s \tilde t)^3=\tilde s^2,\\
 \tilde s^2 \tilde t=\tilde t  \tilde s^2.\endgathered\tag 1.16$$

Since $c^2=1$, this shows that
 $\tilde s, \tilde t$ give a  representation of the modular
group $SL_2(\Z)$. Recall that $SL_2(\Z)$ is generated
by the elements
$$S=\pmatrix 0&-1\\ 1 &0\endpmatrix, T=\pmatrix 1&1\\ 0
&1\endpmatrix$$
satisfying the defining relations $(ST)^3=S^2, S^2T=TS^2, S^4=1$.

Now, let us define the following object in $\Cal C$:

$$H=\bigoplus_{i\in I} X_i\o X_{i^*}.\tag 1.17$$

We assume that we have fixed isomorphisms $X_i^*\simeq X_{i^*}$ as
in (1.1), and thus we can also write $H$  as
$\bigoplus X_i\o X_i^*$ or $\bigoplus X_i^*\o X_i$. Note that since
$(X_i\o X_{i^*})^*\simeq X_{i}\o X_{i^*}$, we have an isomorphism
$H\simeq H^*$.

\definition{Definition 1.8} Define $S, T, C\in \text{End }H$ as follows:
$S=\bigoplus S_{ij}, S_{ij}:X_j\o X_{j^*}\to X_i\o X_{i^*}$,
and similarly for $T, C$, where $S_{ij}, T_{ij}, C_{ij}$ are
given by

$$
S_{ij}=\frac{\dim X_j}{D}\quad \vcenter{\epsfbox{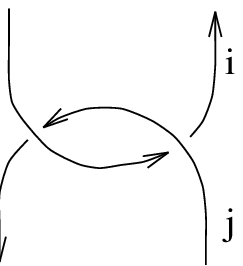}}\qquad
T_{ij}=\frac{\delta_{ij}}{\zeta}\quad\vcenter{\epsfbox{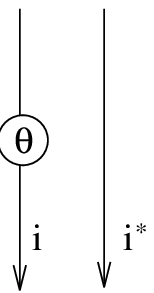}}\qquad
C_{ij}=\delta_{ij^*}\quad\vcenter{\epsfbox{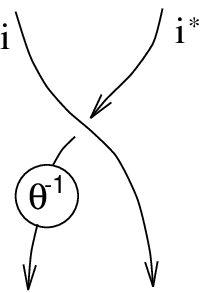}}
\tag 1.18$$

\enddefinition

\proclaim{Theorem 1.9} The morphisms $S,T, C$ defined above satisfy the
following relations:

$$\gathered
S^2= C,\\
S^4=C^2=\theta_H^{-1},\\
(ST)^3= S^2,\\
S^2T=TS^2.\endgathered\tag 1.19$$
\endproclaim
\demo{Proof} Follows from Proposition 1.5 and Corollary~1.6.\enddemo

We cannot say that $S, T$ give a projective representation of
the modular group in $H$, since $\theta_H$ is not a constant. However,
it is true if we restrict them to an isotypic component of
$H$. Equivalently,
let us fix a simple object $U$ in our category and
 consider the space

$$\Hom (H, U)=\bigoplus \Hom (X_i\o X_{i^*}, U).$$

 This is a linear space over $\K$, and $\theta_H |_{\Hom(H, U)}=
\theta_U \Id_{ \Hom(H, U)}$.

\proclaim{Theorem 1.10} Define the maps $S_U, T_U: {\Hom(H, U)}\to
{\Hom(H, U)}$ by

$$ \gathered
S_U: \Phi\mapsto \Phi  S,\\
T_U: \Phi\mapsto \Phi  T.\endgathered\tag 1.20$$

Then $S_U, T_U$ satisfy the following relations

$$\gathered
S_U^4=\theta_U^{-1},\\
T_US_U^2=S_U^2T_U,\\
(S_UT_U)^3=S_U^2,\endgathered\tag 1.21$$
and thus give a projective representation of the  group
$SL_2(\Z)$ in $\Hom (H,U)$. \endproclaim

\example{Example} Let $U=\1$ be the unit object in $\Cal C$. Then we
have a canonical identification $\Hom(X_i^*\o X_i, \1)\simeq \K$, and
thus we have a canonical basis  $\chi_i\in\Hom (H, \1)$. In this case,
the action of the modular group defined in Theorem~1.10 in the basis
$\chi_i$ is given by $\tilde s, \tilde t$
defined by (1.15).
\endexample

The following result, which is a reformulation of theorem of Vafa, is
also worth mentioning here (though we won't use it):

\proclaim{Theorem 1.11} In any modular tensor category \rom{(}regardless of
the base field\rom{)} the numbers $\theta_i, \zeta$ \rom{(}see \rom{(1.14))}
are roots of unity.\endproclaim

This theorem was proved by Vafa (see \cite{Vaf}) in the context of
conformal field theory. However, his proof only uses some relation in
the mapping class group of $n$-punctured sphere and action of
$SL_2(\Z)$. Both of them act in arbitrary modular tensor category: the
action of $SL_2(\Z)$ was discussed above, and the action of mapping
class group can be defined as well (see \cite{Tu, V.4}). Thus, the
same proof is valid in arbitrary MTC. Note that for MTC's coming form
Conformal Field Theory, we have  $\zeta=e^{2\pi\i c/24}$, where $c$ is
the central charge of the action of Virasoro algebra, and theorem
above implies that $c$ is rational, which is why these theories are called
rational.

\subhead Hermitian categories \endsubhead

We will also need the notion of hermitian category: this
definition and all the properties we are citing are due to Turaev (see
\cite{Tu}). Let us
assume that $\Cal C$ is a ribbon category which is defined over
the ground  field $\K$ which is  equipped with an involution
$x\mapsto \bar x$; our basic examples of such an involution will be
$\K=\C$ with usual complex conjugation, and $\K=\C(q), \bar q=
q^{-1}$.
 We say that  $\Cal C$ is {\it  hermitian} if for
every objects $V, W$ we have an involutive
 map $\overline{\phantom{T}}: \Hom (V, W)\to \Hom (V^*,
W^*)$,  such that $\overline{f+g}=\bar f+\bar g,
\overline{\alpha f}=\bar\alpha\bar f$ for any
$\alpha\in \K$,
$\overline{fg}=\bar f\bar g, \overline{f\o g}=\bar g\o \bar f$,
and $\overline{\Id_V}=\Id_{V^*}$. Note
that since we have a canonical identification $\Hom(V^*, W^*)\simeq
\Hom(W, V)$ we could as well consider $\bar f$ as an element of $\Hom
(W, V)$. This involution must satisfy certain compatibility
properties, namely:

$$\gathered
\bar \theta_V=\theta_{V^*}^{-1},\\
\overline{\check R_{V, W}}= \bigl(\check R_{V^*, W^*}\bigr)^{-1},\\
\overline{ e_V}=e_{V}(1\o \delta_V^{-1}):V^*\o V^{**}\to \1,
\endgathered\tag 1.22$$
where $e_V$ is the canonical morphism $V^*\o V\to  \1$.
 Then it can be shown that
geometrically this involution  corresponds to reflection: if $f$ is a morphism
corresponding to the ribbon graph $\Gamma$ then $\bar f$ corresponds
to the graph $\overline\Gamma$ obtained by reflection of $\Gamma$ around a
plane $x=1$ (we assume that the graph is drawn in the projection to
$x,y$-plane) and
changing each label $V$ by $V^*$. Note that this operation changes the
orientation of $\R^3$.

If $\Cal C$ is a hermitian modular category then it follows from the
above geometric interpretation of bar conjugation that we have the
 following identities:

$$ \gathered
\overline{s_{ij}}= s_{i j^*},\\
\overline{\theta_i}=\theta_i^{-1}.\endgathered\tag 1.23$$
 Thus,

$$\gathered
\overline{\dim V}=\dim V,\\
\overline{p^+}=p^-, \quad \overline{p^-}=p^+.\endgathered\tag 1.24$$

Therefore,  $p^+p^-$ is ``real'' and $ p^+/p^-$ is ``unitary'':
$\overline{p^+ p^-}=p^+p^-, \overline{p^+/p^-}=(p^+/p^-)^{-1}$.
We assume that $D$ and $\zeta$ (see (1.14)) can also
be chosen ``real'' and ``unitary'' respectively:

$$\overline{D}=D, \quad\overline{\zeta}=\zeta^{-1}.\tag 1.25$$

Obviously, it is so if $\K=\C$, since in this case $D^2=\sum (\dim
X_i)^2\in \R_+$, and $|\zeta|=1$.

\proclaim{Proposition 1.12} The matrices  $\tilde s, \tilde t\in
Mat_{|I|}(\K)$ are ``unitary'', i.e. satisfy
$X X^*=1$, where $(X^*)_{ij}=\overline{X_{ji}}$.
\endproclaim
\demo{Proof} Obvious from (1.16), (1.23). \enddemo

Similar statement holds in more general case:
\proclaim{Proposition 1.13} The operators $S, T\in \End H$
satisfy the following properties:

$$\gathered
\overline {S}= S C^{-1},\\
\overline{T}= T^{-1}.\endgathered$$

Thus, both $S, T$ satisfy

$$X \overline X =\Id_H.$$
\endproclaim

\demo{Proof} It is obvious for $T$; as for $S$, we need to prove that
$(SC^{-1})_{ij}$
is given by the following picture

$$(SC^{-1})_{ij}=\frac{\dim X_j}{D}\quad\vcenter{\epsfbox{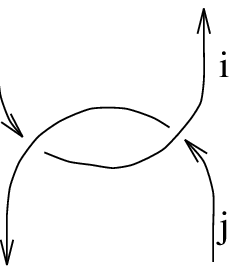}}$$
which easily follows from the definitions. \qed\enddemo

In the next chapter we will show how Proposition~1.13 can be
interpreted as ``unitarity'' with respect to a certain inner product
in $H$.

\head 2. Inner product on morphisms. \endhead

In this section we define an inner product on spaces of morphisms in a
hermitian MTC; this definition is due to Turaev. We also show that the
action of modular group defined in Section~1 is unitary with respect
to this action; the same applies to associativity and commutativity
isomorphisms (when rewritten in terms of $\Hom$ spaces). These results
are new.

As before, we assume that $\Cal C$ is a modular category; we keep all
the notations and conventions of Section~1.

\definition{Definition 2.1} Let $V, W$ be objects from $\Cal
C$. Assume that $\dim V\ne 0, \dim W\ne 0$ and that there exist
$\sqrt{\dim V}, \sqrt{\dim W}$ in $\K$.
 Define the pairing
$$\Hom(V, W)\o \Hom(V^*, W^*)\to \K$$
as follows: if $\Phi_1\in \Hom (V,W)$, $\Phi_2\in \Hom (V^*, W^*)$
then let

$$
\<\Phi_1, \Phi_2\>=\frac{1}{(\dim V \dim
W)^{1/2}}\quad\vcenter{\epsfbox{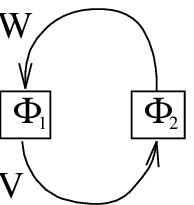}}\tag 2.1$$

\enddefinition
Obviously,  this pairing is symmetric.

\example{Examples}
\roster\item
Let $V=W$. Then $\<\Id_V, \Id_{V^*}\>=1$ (this justifies the choice of
normalization in Definition~2.1).

\item Consider intertwiners of the form $\Phi_1: X_i\o X_j\to X_k,
\Phi_2: X_{j^*}\o X_{i^*}\to X_{k^*}$. Then Definition~2.1 allows to
define pairing between them provided that we have chosen
identifications $X_{i^*}\simeq X_i^*$, etc. Note that in this case
dimension $\dim (X_i\o X_j \o X^*_k)$ is non-zero automatically.
\endroster
\endexample

If $\Cal C$ is a hermitian category then we define a
``hermitian'' inner product on $\Hom(V,W)$ by

$$(\Phi_1, \Phi_2)=\<\Phi_1, \overline{\Phi_2}\>;\tag 2.2$$
as usual, we will denote $\|\Phi\|^2=(\Phi, \Phi)$.
It is easy to see that this inner product satisfies the usual
relations

$$\gathered
(\a\Phi_1, \Phi_2)=\a(\Phi_1, \Phi_2), \quad \a\in \K\\
(\Phi_2, \Phi_1)=\overline{(\Phi_1, \Phi_2)}.\endgathered$$

\proclaim{Lemma 2.2}{\rm (\cite{Tu})}
 In a hermitian modular category, the inner
product  given by \rom{(2.2)} is
non-degenerate. \endproclaim

\remark{Remark} Obviously, the definition of the pairing (and thus, of
the inner product) works as well in a ribbon category without the
assumption of modularity; however, in this case it is not true that
this inner product  is non-degenerate. \endremark

\proclaim{Lemma 2.3} Let $\Phi_1,\Phi'_1: V_1\o V_2\to X_i,
\Phi_2,\Phi'_2:X_i \o
V_3\to U$ be morphisms in a hermitian modular category, and let
$\Psi=\Phi_2(\Phi_1\o 1),\Psi'=\Phi'_2(\Phi'_1\o 1)\in
\Hom(V_1\o V_2\o V_3, U)$. Then

 $$(\Psi,\Psi')=(\Phi_1,\Phi'_1)(\Phi_2,\Phi'_2).\tag 2.3$$
\endproclaim
\demo{Proof} Follows from the identity
$$\vcenter{\epsfbox{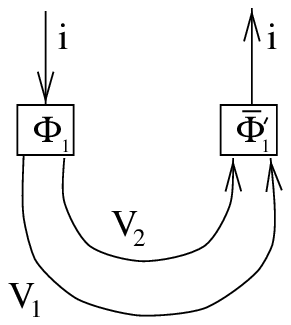}}
=\left(\frac{\dim V_1 \dim V_2}{\dim X_i}\right)^{1/2}
(\Phi_1, \Phi'_1)
\quad\vcenter{\epsfbox{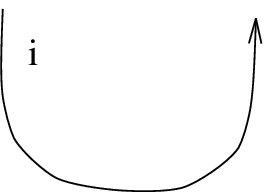}}.$$
\enddemo

We can rewrite commutativity and associativity isomorphisms in terms
of $\Hom$ spaces, which gives us isomorphisms

$$\gathered
\check R:\Hom(V_1\o V_2, U)\to \Hom(V_2\o V_1, U),\\
\a:\bigoplus_{i\in I}\Hom(V_1\o V_2, X_i)\o \Hom(X_i\o V_3, U)\to \\
\qquad \bigoplus_{i\in I}\Hom(V_1\o X_i, U)\o \Hom(V_2\o V_3, X_i).
\endgathered\tag 2.4
$$

\proclaim{Theorem 2.4} In a hermitian modular tensor category, the
associativity and commutativity isomorphisms \rom{(2.4)} are unitary,
i.e. preserve the inner product \rom{(2.2)}. The same is true for the
isomorphism $\Hom(V_1\o V_2, V_3)\simeq \Hom (V_1\o V_2\o {}^*V_3,
\1)$.
\endproclaim

\demo{Proof} Unitarity of associativity isomorphism follows from
Lemma~2.3; unitarity of commutativity follows from the following picture:
\midinsert
\centerline{\epsfbox{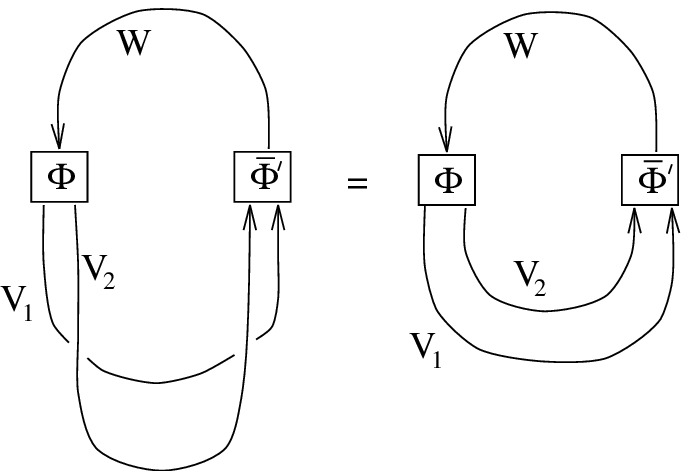}}
\endinsert

Unitarity of the last  isomorphism is obvious.\qed
\enddemo

In particular, (2.2)  gives a natural inner product on
each of the spaces $\Hom (X_i\o X_{i^*}, U)$, and thus -- by taking
direct sum over all $i$ -- on $\Hom (H, U)$. (Note that because of the
normalizations, the inner product $(\Phi_1, \Phi_2)$ depends on
whether we consider $\Phi_1, \Phi_2$ as intertwiners $X_i\o X_{i^*}\to
U$ or  $H\to U$. We will always use the former choice, i.e. consider
$\Phi$ as morphisms $X_i\o X_{i^*}\to U$.)

\proclaim{Theorem 2.5} Let $U$ be an irreducible object. Then the
inner product in the space $\Hom(H, U)$ defined by  \rom{(2.2)} is
invariant under the projective action of the modular group on $\Hom(H,
U)$,  which was defined in Theorem~\rom{1.10}.
\endproclaim

\demo{Proof} In view of the identities $S_U \overline{S_U}=1, T_U
\overline{T_U} =1$ (Proposition 1.13), it  suffices to show that
$(\Phi_1 S, \Phi_2)=(\Phi_1, \Phi_2\overline{ S})$, or,
equivalently, $\<\Phi_1 S, \overline{\Phi_2}\>=\<\Phi_1,
\overline{\Phi_2}S\>$, and
similarly for $T$. For $T$ this is obvious; for $S$, it follows  from
the  following picture:

$$\epsfbox{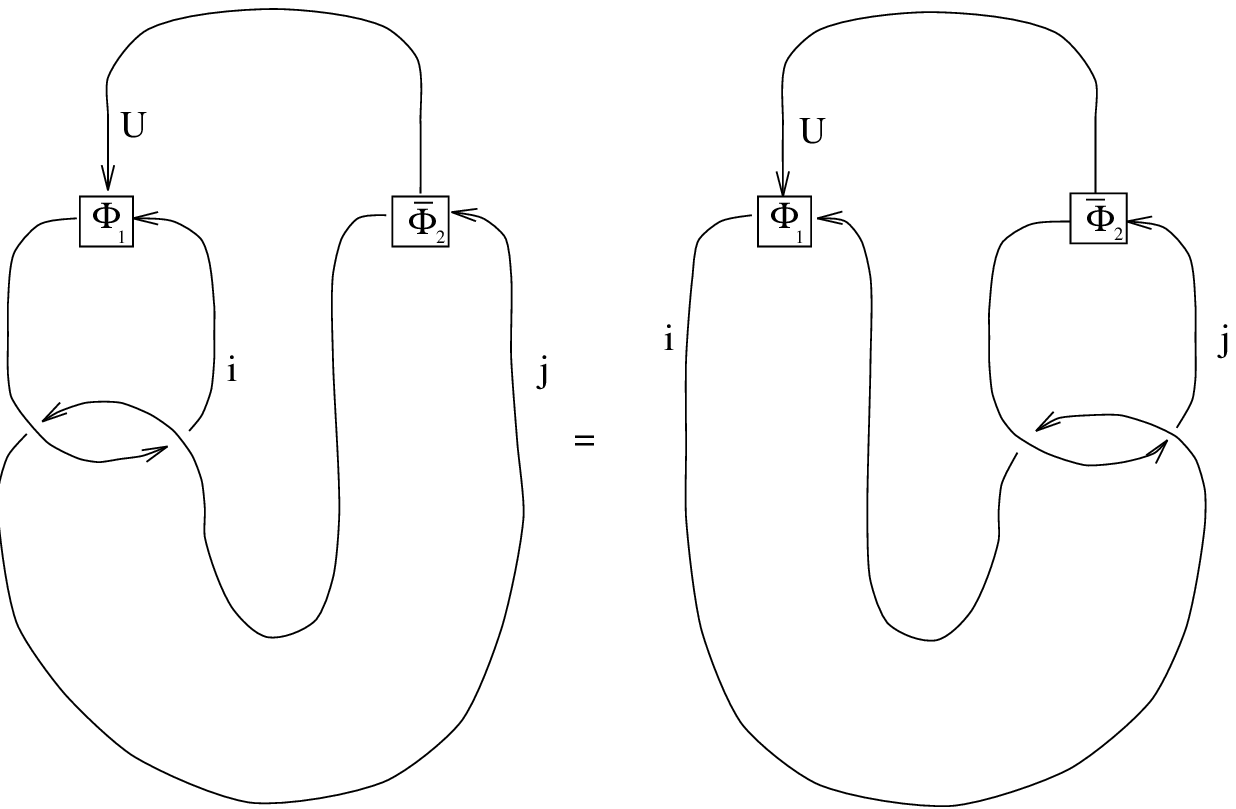}$$
 \enddemo

\head 3. Quantum groups at roots of unity\endhead

In this section we recall the known results on construction of modular
categories from representations of quantum groups at roots of unity,
following the papers of Andersen
(\cite{A, AP}) (which, in turn, are based on the work of Lusztig
 (\cite{L1--L4}, see also \cite{L5})), and Gelfand-Kazhdan
(\cite{GK}). Again, this section is completely expository.

\subhead General facts on quantum groups\endsubhead

Here we give the main definitions from the theory of quantum groups;
they are well known, and we refer the reader to the original papers by
Drinfeld and Jimbo or to Lusztig's book \cite{L5} for details, giving here the
bare minimum -- mostly to fix notations.

Let $\g$ be a simple Lie algebra over $\C$. We use the standard
notations for the Cartan
subalgebra, roots, weights etc.; we also denote by $\theta$ the
highest root of $\g$. We denote by $(\, , \,)$ invariant
bilinear form in $\h$ normalized so that $(\theta, \theta)=2$ and by $(\, ,
\,)'=m(\, , \,)$ the form normalized so that $(\a, \a)'=2$ for short
roots; equivalently, it is specified by the conditions $d_i=
(\a_i, \a_i)'/2\in \Z_+, \gcd d_i=1$. Thus, $m=1$ for simply-laced
root systems, $m=2$ for root systems of $B, C, F$ types and $m=3$ for
$G_2$.

 By definition,
the corresponding quantum group  $\Ug$  is an algebra over the field
$\C_q= \C(q^{1/2N})$, where $N=|P/Q|$ (fractional powers are necessary
to define braiding)
with generators $e_i, f_i, q^{h},h\in \frac{1}{2}Q\v\subset \h,
 i=1\ldots r$ and relations

$$\gathered
q^0=1, \quad q^{a+b} =q^{a}q^b,\\
q^h f_j q^{-h}=q^{-\<h, \a_j\>} f_j,\\
q^h e_j q^{-h}=q^{\<h, \a_j\>} e_j,\\
[e_i, f_j]=\delta_{ij}\frac{q^{d_i h_i}-q^{-d_i h_i}}{q^{d_i}-q^{-d_i}},
\endgathered \tag 3.1$$

$$\gathered
\sum_{n=0}^{1-a_{ij}}
	\frac{(-1)^n}{[n]_i! [1-a_{ij}-n]_i!}
	e_i^n e_j e_i^{1-a_{ij}-n}=0,\\
\sum_{n=0}^{1-a_{ij}}
	\frac{(-1)^n}{[n]_i! [1-a_{ij}-n]_i!}
	f_i^n f_j f_i^{1-a_{ij}-n}=0,\endgathered
\tag 3.2$$

where  $h_i=\a_i\v\in \h$ and
$$[n]_i=\frac{q^{nd_i} - q^{-nd_i}}{q^{d_i}-q^{-d_i}},
\quad [n]_i! =[1]_i \dots [n]_i.\tag 3.3$$

This is a Hopf algebra with the following comultiplication, counit and
antipode:
$$\gathered
\Delta e_i =e_i\o q^{d_i h_i/2} +q^{-d_i h_i/2}\o e_i,\
\Delta f_i =f_i\o q^{d_i h_i/2} +q^{-d_i h_i/2}\o f_i,\\
\Delta q^h=q^h\o q^h,\\
\epsilon(q^{h})=1, \epsilon(e_i)=\epsilon(f_i)=0,
S(e_i)=-q^{d_i}e_i, S(f_i)=-q^{-d_i}f_i, S(q^h)=q^{-h}.
\endgathered\tag 3.4$$

As is well-known, this algebra is quasitriangular: there exists a
``universal R-matrix'' $\Cal R$ which is an element of a certain completion of
$ \Ug\o \Ug$
such that for every pair of finite-dimensional representations $V,W$
the operator
$$\check R_{V,W}= P\circ \pi_V\o \pi_W(\Cal R) \colon V\o W\to
W\o V
\tag 3.5$$
is an isomorphism of representations. Here $P$ is the transposition:
$Pv\o w=w\o v$.
Also, it is known that $\Cal R$ has the following
form:
$$\gather
\Cal R=q^{\sum a_i\o a_i}\Cal R^*, \quad
	\Cal R^*\in U^+\hat \o  U^-\tag 3.6\\
(\epsilon\o 1)(\Cal R^*)=(1\o\epsilon)(\Cal R^*)=1\o
1,\endgather$$
where $a_i$ is an orthonormal basis in $\h$ with respect to $(\, ,\,)'$.
As we said, $\Cal R$ does not lie in the tensor square of
$\Ug$ but in its certain completion; however, for any pair
of finite-dimensional representations $V, W$ of $\Ug$ the
operator $\pi_V \o \pi_W(\Cal R)$ is well-defined (this is
where we need fractional powers of $q$ in the definition of
$\C_q$).

\remark{Remark} This definition differs from Lusztig's one by slightly
different choice of generators and, more importantly, by replacing
$v=q^{-1}$. \endremark

We recall (see \cite{Kas, Tu})
that the category $Rep \ \Ug$ of finite-dimensional
representations of $\Ug$ is a semisimple ribbon category (in the sense
of definitions of Section~1) defined over $\C_q$. Its simple objects
are precisely irreducible highest-weight modules
$V_\l, \l\in P^+$.
Note that if  $\l\in P^+$ then
$V_\l^*\simeq V_{\l^*}$, where $\l^*=-w_0(\l)$, $w_0$ being the
longest element of the Weyl group $W$. In this category,  the balancing map
 is such that  $\theta_{V_\l}=q^{(\l, \l+2\rho)'}$, and the isomorphism
$\delta_V: V\to V^{**}$ is given by $q^{2\rho}$, where $\rho$ is
considered as an element of $\h$ using the identification given by
$(\, , \,)'$; thus, $q^{2\rho}v = q^{2(\rho, \l)'}v$ if $v$ has weight
$\l$. This implies that the quantum dimension of a module is
given by $\dim_q V= \Tr_V (q^{2\rho})$. In particular, if $V_\l$ is
irreducible then $\dim_q V_\l=\chi_\l(q^{2\rho})$, where $\chi_\l\in
\C[P]$ is the character of $V_\l$, and we use the following
convention: for
$f=\sum a_\l e^\l\in \C_q[P]$ we let

$$f(q^\mu)=\sum a_\l q^{(\l, \mu)'}.\tag 3.7 $$

It follows from Weyl formula that
$$\dim_q V_\l=\chi_\l(q^{2\rho})=\frac{\sum_{w\in W}
(-1)^{l(w)}q^{2(\rho,w(\l+\rho))'}}{\delta(q^{2\rho})}
=\frac{\delta(q^{2(\l+\rho)})}{\delta(q^{2\rho})},
\tag 3.8$$
where $\delta$ is Weyl denominator:

$$\delta=\sum_W (-1)^{l(w)} e^{w(\rho)}=\prod_{\a\in
R^+}(e^{\a/2}-e^{-\a/2}).\tag 3.9$$

\subhead Representations of $\Ug$ at roots of unity and category $\Cal
C(\g, \varkappa)$\endsubhead

Let $A=\Z[q^{1/2N}, q^{-1/2N}]$, and let $U$ be the
$A$-subalgebra of $\Ug$ generated by $e_i^n/[n]!$, $f_i^n/[n]!, q^h$
 (see \cite{L2}). For arbitrary non-zero number $\eps\in \C$ define

$$U_\eps=U\o_A \C,\tag 3.10$$
where $\C$ is endowed with a structure of an $A$-module by $q\mapsto
\eps$.

Our goal is to construct a certain subquotient of the category of
finite-dimensional representations of $U_\eps$ in the case when $\eps$
is a root of unity:

$$\eps=e^{\pi \i /m\varkappa}, \tag 3.11$$
where $m$ is as before and $\varkappa\in \Z_+$.
In this paper we always
assume that $\varkappa\ge h\v$, where $h\v$ is the dual Coxeter number
for $\g$. By definition, we  let $\eps^a= e^{\pi\i a/m\varkappa}$ for
$a\in \frac{1}{2N}\Z$, and as before,
for
$f=\sum a_\l e^\l\in \C[P]$ we let

$$f(\eps^\mu)=\sum a_\l \eps^{(\l, \mu)'}= \sum a_\l e^{\pi \i (\l,
\mu)/\varkappa}.\tag 3.12$$

Recall (see \cite{L1}) that due to the fact that $U_\eps$ contains
divided powers, we have a notion of weight subspace, and weight
subspaces are indexed by $P$ (not by $P/2\varkappa P$!). Let
$ Rep\ U_\eps$ be the category of finite-dimensional representations of
$U_\eps$ with weight decomposition.

\proclaim{Theorem  3.1}  $ Rep\
U_\eps$ has a natural structure of ribbon category over $\C$. \endproclaim

\demo{Proof} This follows from general result due to Lusztig
(\cite{L5, Chapter 32}); we refer the reader to \cite{KL4, \S37} for
details.\enddemo

 Let $V_\l$ be the
irreducible finite-dimensional  module over $\Ug$ with highest weight
$\l\in P^+$, and let $v_\l$ be the highest weight vector in it. $V_\l$
admits a $U$-structure: we can consider $U$-submodule $U v_\l\subset
V_\l$. Thus, we can define  a module over $U_\eps$ which we
denote by $V^\eps_\l$ (sometimes we will also denote it by $V_\l$):

$$V^\eps_\l=(U v_\l)\o_A \C.\tag 3.13$$

These modules are usually called Weyl modules and are not necessarily
irreducible (see below).

Define the open and closed alcoves $C, \overline{C}$ by
$$ \gathered
C=\{\l\in P^+| \<\l+\rho,
\theta\v\> < \varkappa\},\\
\overline{C}=\{\l\in P|\<\l+\rho, \a_i\v\>\ge 0,
\<\l+\rho, \theta\v\>\le\varkappa\}.\endgathered
 \tag 3.14$$
Note  that $C$ is preserved by the involution $\l\mapsto \l^*=-w_0(\l)$.

We will denote by $\Gamma$ be the affine wall of $\overline{C}$:

$$\Gamma=\{\l\in\h^*|\<\l+\rho, \theta\v\>=\varkappa\}.\tag 3.15$$

Note that $\overline{C}$ is the fundamental domain for the shifted action of
affine Weyl group $\What$ of level $\varkappa$. Recall that the shifted
action is defined by  $w.\l=w(\l+\rho)-\rho$, and
 $\What=W\ltimes
\varkappa Q\v$, where $Q\v$ is considered as a lattice in $\h^*$ using
the identification  $\h\simeq \h^*$ given by  the form $(\, , \,)$;
under this identification $Q\v\subset Q$.

Now we can formulate the main result on  the reducibility of Weyl modules:

\proclaim{Lemma 3.2} For $\l\in
\overline{C}\cap P^+$, Weyl modules $V_\l$ are irreducible.

\endproclaim

In general,  Weyl modules are not irreducible.
However, it is easy to see that for every $\l\in P^+$ there exists a
unique irreducible highest-weight module $L_\l$, and that every
irreducible module in $Rep\ U_\eps$ is of the form $L_\l, \l\in P^+$.
Thus, every module from $Rep \ U_\eps$ has a composition series with
factors of the form $L_\l$. In particular, the same is true for $V_\l$,
and the multiplicities  should (conjecturally) express in terms of
Kazhdan--Lusztig polynomials associated with the affine Weyl group
$\What$ (see \cite{L6, Section~9}).

Our goal is to extract from this highly non-trivial category of
representations some semisimple part. This was first done by
Reshetikhin and Turaev in the case $\g=\sltwo$ and by Andersen
(\cite{A, AP}) in general case.  We
briefly sketch the main  steps here.

Let us call a module $V\in Rep \ U_\eps$ {\it tilting} if both $V$ and
$V^*$ have  a composition series with
the factors isomorphic to $V_\l$. Let $\Cal T$ be the full subcategory
of $ Rep \  U_\eps$ consisting of tilting modules. This category
 is closed under taking
dual representations (obvious) and under tensor product (see
\cite{A, AP}), and thus is a ribbon category. Note that in particular
the modules $V_\l, \l\in C$ are tilting.

However, the category of tilting modules is still too large, and thus
we need further reduce it. This is done by factorization over
negligible modules.

For every finite-dimensional module $V$ over $U_\eps$ we define its
dimension $\dim_\eps V$ by the same formula as for $\Ug$; in
particular, for the highest-weight module $V^\eps_\l$ we have $\dim_\eps
 V^\eps_\l= \chi_\l(\eps^{2\rho})$. Simple calculation gives the
following result:

\proclaim{Lemma 3.3} Let $\l\in P^+$. Then

$$\dim_\eps V_\l =0\iff
(\l+\rho, \a)\in \varkappa\Z\text{ for some }\a\in R^+.$$

 In particular,
if $\l\in C$ then $\dim_\eps V_\l\ne 0$, and if $\l\in \Gamma\cap P^+$
then $\dim_\eps V_\l=0$ \rom{(}recall that $\Gamma$ denotes the affine
wall of $\overline{C}$\rom{)}.\endproclaim

Let us call a tilting module $V$ {\it negligible} if for every
$f\in \End V$ we have $\Tr_\eps(f)=0$. The following simple lemma
describes the properties of negligible modules:

\proclaim{Lemma 3.4}{\rm(\cite{A})}
\roster\item
 An indecomposable  module $V\in Rep\ U_\eps$ is
negligible iff $\dim_\eps V=0$.

\item $V$ is negligible iff $V^*$ is negligible.

\item If $\l\in C$ then $V_\l$ is not a direct summand of a negligible
module; equivalently, for every $\l\in C$ and a negligible
module $Z$ the composition

$$V_\l @>{f}>> Z@>{g}>>V_\l$$
is equal to zero for any morphisms $f,g$.
\item If $V$ is a negligible module then $V\o V'$ is negligible for
any $V'\in Rep \ U_\eps$.
\endroster
\endproclaim

The following key theorem is due to Andersen:

\proclaim{Theorem 3.5} Every tilting  module $V$ can be written in
the form

$$V=\left(\bigoplus_{\l\in C} n_\l V_\l\right) \oplus Z$$
for some negligible tilting module $Z$ and uniquely defined non-negative
integers $n_\l$.
\endproclaim

In particular, this theorem implies that for $\l, \mu \in C$
we can write

$$V_\l \o V_\mu=\biggl(\bigoplus_{\nu\in C} N_{\l\mu}^\nu V_\nu\biggr)\oplus
Z,\tag 3.16$$
for some $N_{\l\mu}^\nu\in \Z_+$ and $Z$ as above.

In the case $\g=\sltwo$ this theorem was proved by  Reshetikhin and
Turaev.

This theorem allows us to define the modular category with simple
objects $V_\l, \l \in C$ as follows.
Recall that $\Cal T$ is the full subcategory of tilting modules in
$Rep \  U_\eps$. Let
$\Cal T^{neg}\subset \Cal T$ be the full  subcategory of
negligible tilting modules. We want to  define the quotient category
$\Cal C=\Cal T/\Cal T^{neg}$; this construction is due to
\cite{GK}, and we briefly sketch it below.

Let $V_1, V_2\in \Cal T$. We call a
morphism $f:V_1\to V_2$ negligible if it can be presented in the form
 $f=gh$ for some $h:V_1\to Z, g:Z\to V_2$, where
$Z$ is a negligible module. We
denote negligible morphisms from $V$ to $W$ by $\Hom^{neg}(V,W)$.

\definition{Definition 3.6} The quotient category $\Cal C(\g,
\varkappa)=\Cal C$
is defined as follows:
$$\gathered
\text{\rm Ob } \Cal C=\text{\rm Ob }  \Cal T, \\
\Hom_{\Cal C}(V, W)=\Hom_{ U_\eps}(V,
W)/\Hom^{neg}_ {U_\eps}(V,W).\endgathered $$
\enddefinition

It follows from Lemma~3.4 that
if $f$ is a negligible morphism then for any morphism $g$ the
composition $fg$ is also negligible; the same applies to $gf, f^*$ and
$f\o g$. Therefore, compositions, tensor products and duals  of
morphisms   are  well-defined and thus  $\Cal C(\g,
\varkappa)$ is a ribbon category. Obviously, in this
category   every negligible module is isomorphic to the  zero module. Thus,
Theorem~3.5 implies
that every object in $\Cal C$ is isomorphic to a direct sum

$$V=\bigoplus_{\l\in C} n_\l V_\l$$
for some unique collection of non-negative integers $n_\l$. In
particular, we have the following isomorphism in $\Cal C$:

$$V_\l\o V_\mu\simeq \bigoplus_{\nu\in C} N_{\l\mu}^\nu V_\nu, \tag 3.17$$
where the numbers $N_{\l\mu}^\nu$ are the same as in (3.16).

{\bf Warning:} in general, $N_{\l\mu}^\nu\ne \dim \Hom_{U_\eps}
(V_\l\o V_\mu, V_\nu)$. Instead, we have the following result:

\proclaim{Lemma 3.7} Let $\l,\mu,\nu\in C$. Then

$$\Hom_{U_q \g}(V_\l\o V_\mu, V_\nu)
\text{\rm `` = ''}\Hom_{U_\eps}(V_\l^\eps\o
V_\mu^\eps, V_\nu^\eps).$$

The equality should be understood in the following sense: we can
define some intertwining operators $\Phi_i$ which are defined over
$A=\Z[q^{\pm 1/2N}]$ such that $\Phi_i$, considered as intertwining
operators over $\C_q$ \rom{(}respectively, over $\C$\rom{)} form a basis in
$\Hom_{U_q \g}(V_\l\o V_\mu, V_\nu)$ \rom{(}respectively, in
$\Hom_{U_\eps}(V_\l^\eps\o V_\mu^\eps, V_\nu^\eps)$\rom{)} -- compare
with definition \rom{(3.13)} of $V_\l^\eps$.
\endproclaim

\demo{Proof} It follows from the fact that one can write explicit
formula for such an intertwiner involving only the inverse of the
Shapovalov form (see, for example, \cite{EK3}).
Since Shapovalov form is non-degenerate in both
$V_\l$ (as a matrix with entries from $\C_q$) and in $V_\l^\eps$ (as a
matrix with complex entries), this proves the lemma.\qed \enddemo

\subhead $\Cal C(\g, \varkappa)$ as a modular category\endsubhead

Let us summarize the properties of the category $\Cal C$:

\proclaim{Proposition 3.8}
\roster\item The category $\Cal C(\g, \varkappa)$ is semisimple, and
its simple objects are precisely $\{V_\l\}_{\l\in C}$.

\item For any object $V\in \Cal C(\g, \varkappa)$ we have $\dim_\eps
V\in \R_{>0}$.

\item This category has a natural structure of ribbon category,
 inherited from $Rep \ U_\eps$.

\item The matrices $s_{ij}, t_{ij}$ defined by \rom{(1.4), (1.12)} for
the category  $\Cal C(\g, \varkappa)$ are  given by

$$\gathered
t_{\l\mu}=\delta_{\l\mu}\eps^{(\l, \l+2\rho)'},\\
s_{\l\mu}=\chi_\l(\eps^{-2(\mu+\rho)})\dim_\eps V_\mu^\eps
=\frac{\sum_{w\in W}(-1)^{l(w)}\eps^{-2(w(\l+\rho), \mu+\rho)}}
	{\delta(\eps^{-2\rho})},\endgathered
\tag 3.18$$
where $\chi_\l\in \C[P]^W$ is the  character of the module $V_\l$,
 $\delta$  is Weyl denominator  and we use convention \rom{(3.12)}.
\endroster
\endproclaim

\demo{Proof} (1) was already discussed; (2) follows from  Weyl
formula for $\dim_\eps V_\l$ and (1); (3) is obvious. Formula (3.18)
is also very well-known and can be deduced from the diagonal part of
the $R$-matrix (see, for example, arguments in \cite{EK3}).
\qed\enddemo

\proclaim{Theorem 3.9} The matrix $s$ defined by \rom{(3.18)} is
non-degenerate, and thus $\Cal C(\g, \varkappa)$ is a modular
category. Also, in this case the numbers $D=\sqrt{p^+p^-}
=\sqrt{\sum (\dim_\eps V_\l)^2}$, $\zeta=(p^+/p^-)^{1/6}$
\rom{(}cf. \rom{(1.14))}  are given by

$$\gathered
D=\frac{\sqrt{|P/\varkappa Q\v|}}{\i^{|R^+|}\delta(\eps^{-2\rho})}=
\frac{\sqrt{|P/\varkappa Q\v|}}
	{\prod_{\a\in R^+} 2 \sin \frac{(\a, \rho)}{\varkappa}
\pi},\\
\zeta=\eps^{\frac{\varkappa-h\v}{h\v}(\rho, \rho)'}
= e^{2\pi\i c/24},\quad c= \frac{(\varkappa-h\v)\dim
\g}{\varkappa}. \endgathered
\tag 3.19$$
Here $Q\v$ is considered as a
sublattice in $Q$ via the identification $\h\simeq \h^*$ given by $(\, ,
\,)$.
\endproclaim
\demo{Proof} This follows from the results of Kac and Peterson (see
below) and the ``strange formula of Freudental-de Vries'' (see
\cite{Kac, 12.1.8}):

$$\frac{\dim \g}{24}=\frac{(\rho, \rho)}{2h\v}.$$

 We will give an elementary proof of the identity  $s^2= D^2
c$, where $D$  is given by (3.19) in Section~6. The formula for
$\zeta$ can be proved similarly.
\qed
\enddemo

Note that (3.19) implies the following formula for the renormalized
$s, t$ matrices:

$$\gathered
\tilde s_{\l\mu}=\i^{|R^+|} |P/\varkappa Q\v|^{-1/2}
\sum_{w\in W} (-1)^{l(w)} e^{-2\pi\i (w(\l+\rho),
\mu+\rho)/\varkappa},\\
\tilde t_{\l\mu}=e^{2\pi\i\left(\frac{(\l,
\l+2\rho)}{2\varkappa}-\frac{c}{24}\right)},\endgathered
\tag 3.20$$
where $c$ is given by (3.19).

The  same formulas have appeared as the matrices  of modular transformations
of the characters of integrable modules over affine Lie algebra of level
$k=\varkappa-h\v$ (Kac-Peterson formula, see \cite{Kac,
Section~13.8}). In this case the number $c$ is interpreted as the
central charge of the Virasoro algebra.
We will discuss  the relation between affine Lie algebras and quantum
groups in forthcoming papers.

\remark{Remark 3.10} In fact, the results of this section
are valid in more general
case. Namely, assume that $\varkappa$ is a rational number:
$|\varkappa|=p/q, p, q\ge 1, (p,q)=1$. Then  all
the results above except for Theorem~3.9 are valid with appropriate
changes indicated below. However, if $q$ is not relatively prime with
$|P/Q\v|$ then the matrix $s$ may be degenerate and thus the $\Cal C$ is
not a modular category; however, it is still a semisimple ribbon
category with finite number of simple objects.

\roster\item Assume that $(m,q)=1$. Then the alcove $C$
must be taken to be

$$C=\{\l\in P^+|\<\l+\rho, \theta\v\><p\}=\{\l\in P^+|(\l+\rho,
\a)<p\text{ for all }\a\in R^+\}, \tag 3.21
$$
and we must consider the action of  affine Weyl group
of level $p$ rather than $\varkappa$: $\What=W\ltimes p Q\v$.

\item  Assume that $q$ is divisible by $m$. In this case
we must take
$$C=\{\l\in P^+|\<\l+\rho, \a\v\><p\text{ for all }\a\in R^+\}=
\{\l\in P^+|\<\l+\rho, \phi\v\><p\},\tag 3.22$$
where $\phi\in R^+$ is such that $\phi\v$ is the highest root of
$R\v$, and Weyl group $\What $ should be replaced by the Weyl group
$\What^\natural = W\ltimes p Q$, which is the affine Weyl group
corresponding to $R\v$. In this case the order of $\eps$
 is relatively prime to $m$. This case was  considered
in earlier papers of  Andersen et al. and for prime $p$ it is related
with representations of algebraic groups in characteristic $p$.
\endroster\endremark

\head 4. Hermitian structure on $\Cal C(\g, \varkappa)$. \endhead

In this section we define a hermitian structure on the category $\Cal
C(\g, \varkappa)$ in the sense of Section~1. To the best of my
knowledge, these results are new; however, they are closely related
with the results of Wenzl (see \cite{We}) who considered unitarity of
corresponding representations of Hecke and Birman-Wenzl algebras.

This hermitian structure does not rely on the fact that $q$ is a
root of unity. Therefore, in this section we consider more general
case: $\Ug$ is considered as an algebra
over $\C_q$ with the conjugation $\overline{\phantom{T}}$ in $\C_q$
which extends complex conjugation on $\C$ by  $\overline{q^a}= q^{-a},
a\in \frac{1}{2N}\Z$.

 As we will show, this hermitian structure  is
essentially equivalent to defining an invariant hermitian form on
representations of $\Ug$ satisfying certain conditions. To do it,
we first need to define a structure of $*$-algebra (that is,
a certain involution) on $\Ug$.

Recall that the involution $\l\mapsto -w_0(\l), w_0$ -- the longest
element of the Weyl group, preserves the set of simple roots. Thus, we
have an involution $\vee:[1,\dots, r]\to [1,\dots, r]$ such that
$\a_{i\v}= -w_0(\a_i)$.

\proclaim{Lemma 4.1} There exists  a unique antilinear algebra
automorphism $\omega: \Ug\to \Ug$ such that
$$\aligned
	\omega:         & e_i\mapsto e_{i\v},\\
		       & f_i\mapsto f_{i\v},\\
		       &q^h\mapsto q^{w_0(h)},\\
	               & q\mapsto q^{-1}.
\endaligned\tag 4.1$$

So defined $\omega$ is coalgebra antiautomorphism and satisfies

$$\gathered
\omega^2=1,\quad S\omega=\omega S^{-1},\\
\omega (\Cal R)= \Cal R^{-1}.\endgathered$$
 where  $\omega$ is
extended to $\Ug^{\o 2}$ by $\omega(a\o b)=\omega( a)\o \omega( b)$.
\endproclaim

\demo{Proof} Let $\overline{\phantom{T}}:\Ug\to \Ug$ be the antilinear
involution such that
$\overline{e_i}=e_i,\overline{f_i}=f_i,\overline{q^h}=q^{-h}$ (this is
slightly different from the definition of bar involution in \cite{L5}
due to different choice of generators). One easily checks that this is
a coalgebra antiautomorphism, satisfying $\overline{Sx}=S^{-1}\bar x$ and
$\overline{\Cal R}=\Cal R^{-1}$. Composing it with ($\C_q$-linear)
involution $e_i\mapsto e_{i\v}, f_i\mapsto f_{i\v}, q^h\mapsto
q^{-w_0(h)}$, which obviously preserves all structures of $\Ug$, we
get $\omega$. \qed
\enddemo

Now, for every module $V$ over $\Ug$ define the new module $V^\omega$
as follows: as a set (and more over, as an $\R$-vector space)
$V^\omega$ coincides with $V$, and the action of $\Ug$ is defined by
$\pi_{V^{\omega}}(x)=\pi(\omega x)$. It is easy to see that if
$V=V_\l, \l\in P^+$ then $V^\omega\simeq V_{\l^*}$.
 For a vector $v\in V$ we will
write $v^\omega$ to denote the same vector considered as an element of
$V^\omega$; similarly, if  $\Phi\in \Hom_{\Ug} (V, W)$
 then $\Phi$ is also an intertwiner considered as a
map $V^\omega\to W^\omega$; we will denote it by $\Phi^\omega$.

It follows from the fact that $\omega$ is antiautomorphism of
coalgebras that the map $(v\o w)^\omega\mapsto w^\omega\o v^\omega$ is
an isomorphism $(V\o W)^\omega\simeq W^\omega\o V^\omega$.

In particular, this implies that if $\check R_{V,W}:V\o W\to W\o V$ is
the commutativity isomorphism defined above then

$$\check R_{V,W}^\omega = \bigl(\check R_{V^\omega,W^\omega}\bigr)^{-1}:
W^\omega\o V^\omega\to V^\omega\o W^\omega.\tag 4.2$$

Now, since $V_\l^\omega\simeq V_{\l^*}$, we can identify
$V_\l^\omega\simeq (V_\l)^*$. In other words, there is a unique up to
a constant $\Ug$-homomorphism $V^\omega_\l\o V_\l\to \C_q$, or
a non-degenerate hermitian form $H$ in $V_\l$ such that

$$H(xv,v')=H(v, x^*v'),$$
where $x^*=S\omega(x)$. This form satisfies the usual symmetry condition
$H(v,w)=\overline{H(w,v)}$.

As we said above, this form is defined uniquely up to a non-zero
complex factor, and there is no canonical choice of this form.
Note,  however, that this form can not be positively definite:
if $v\in V[\l], v'\in V[\l']$ then $H(v, v')=0$ unless $\l=w_0(\l')$;
in particular, $H(v,v)=0$ unless $v\in V[0]$.

Since every module is completely reducible, we can choose an
identification $V^\omega\simeq V^*$ for every $V$. Moreover, we can do
it in such a way that this is compatible with tensor product and duality, i.e.
the identification $(V\o W)^\omega\simeq (V\o W)^*$ coincides with the
composition $(V\o W)^\omega\simeq W^\omega\o V^\omega\simeq  W^*\o
V^*\simeq (V\o W)^*$, and identification
$V=V^{\omega\omega}\simeq (V^*)^{\omega}\simeq V^{**}$ coincides with
$\delta_V$.

Thus, if $\Phi$ is an intertwiner $V\to W$ then $\Phi^\omega$ can also
be considered as an intertwiner $V^*\to W^*$, which gives us the
following result:

\proclaim{Theorem 4.2} The map $\Phi\mapsto \Phi^\omega$ defined above
endows  $Rep\ \Ug$
with a structure of hermitian category over the field $\C_q$ with
respect to the above defined complex conjugation on $\C_q$.
\endproclaim
\demo{Proof} We have to check consistency relations (1.22). It follows
from (4.2) that  $\overline{\check R_{V, W}}= \bigl(\check R_{V^*,
W^*}\bigr)^{-1}$; the relation
$\overline{\theta}=\theta^{-1}$ is obvious since $\theta$ has
eigenvalues $q^{(\l, \l+2\rho)'}$, and the commutation relation with
$e_V$ follows from compatibility with duality (see above). \qed\enddemo

The conjugation $\omega$ works as well if we replace
$q$ by a root of unity $\eps$. Moreover,  it is easy to see that
$\omega$ preserves Weyl (and thus, tilting) modules
 and that a morphism $\Phi$ is
negligible iff  $\Phi^\omega$ is negligible. Obviously,
$V_\l^\omega\simeq V_\l^*$ if $\l\in C$; thus, the
construction above defines a structure of hermitian category on $\Cal
C(\g, \varkappa)$.

Having defined the hermitian structure, we can define inner product on
intertwiners $\Hom_{\Cal C(\g, \varkappa)}(V,W)$. In fact, the
construction above
gives even more: it gives an inner product on a larger space
$\Hom_{U_\eps}(V,W)$ if $V, W$ are modules over $U_\eps$.
Recall that by definition we have
$$\Hom_{\Cal C(\g, \varkappa)}(V, W)=
\Hom_{U_\eps}(V, W)/
\Hom^{neg}_{U_\eps}(V, W),$$
where $\Hom^{neg}$ is the space of negligible morphisms.

\proclaim{Lemma 4.3} Let $\l, \mu, \nu\in C$. Then
$\Phi\in \Hom_{U_\eps}(V_\l^\eps\o
V_\mu^\eps, V_\nu^\eps)$ is negligible iff $\Phi$  is in the kernel of
the inner product $(\, , \, )$  defined by \rom{(2.2)}. \endproclaim

\demo{Proof} If $\Phi$ is negligible, then $(\Phi, \Phi')$ can be
rewritten as a trace of some operator in a  negligible module $Z$, and
thus is equal to zero. Vice versa, assume that   $\Phi$ lies in the kernel of
this inner product. Since the inner product of intertwiners in $\Cal
C(\g, \varkappa)$ is non-degenerate (Theorem~2.2), this shows that
$\Phi=0$ as a homomorphism in $\Cal C(\g, \varkappa)$, and thus, by
definition, $\Phi$ is negligible. \qed\enddemo

Finally, let us consider the inner product on the spaces $\Hom (V, V\o
U)$. It turns out that in this case the inner product on intertwiners
coincides with the inner product on  so-called generalized
characters (see \cite{EK1--EK3}), definition of which we briefly
recall below. The arguments below work for both categories $\Cal C(\g,
\varkappa)$ (over $\C$) and for $Rep \ \Ug$ (over $\C_q$); for
simplicity, we will formulate all results for  $\Cal C(\g,
\varkappa)$.

For an intertwiner $\Phi \in \Hom_{U_\eps}(V, V\o U)$ define the
corresponding generalized character $\chi_\Phi\in\C[P]\o U[0]$ by

$$ \chi_\Phi= \sum_{\l\in P} \Tr_{V[\l]}(\Phi).\tag 4.3$$

Equivalently, we can consider $\chi_\Phi$  as a function on $\h$ by
letting $e^\l(h)=e^{\<\l, h\>}$; then the above definition is
equivalent to
  $$ \chi_\Phi(h)=\Tr_V(\Phi e^h). $$

Let us define the following involution on $\C[P]$:

$$\overline{\sum a_\l e^\l}=\sum \overline{a_\l}e^{-w_0(\l)}.\tag
4.4$$

Then one can define the following inner product on $\C[P]\o U[0]$:

$$(\chi_1, \chi_2)_1=\frac{(-1)^{|R^+|}}{|W|}[(\chi_1\o \bar \chi_2)_U
\delta\bar \delta]_0\tag 4.5$$
(the subscript $1$ will be explained later when we generalize this inner
product introducing $(\, , \,)_k$).
Here $\delta$ is Weyl denominator \rom{(3.9)},
 $(\cdot\o \cdot)_U:  (\C[P]\o U[0])^{\o 2}\to
\C[P]$ is composition of the hermitian form $H: U\o U\to \C$
discussed above and multiplication in $\C[P]$, and

$$\left[\sum a_\l e^\l\right]_0=a_0.\tag 4.6$$

Then we have the following theorem which is the hermitian analogue of
statement proved in \cite{EK2}:

\proclaim{Theorem 4.4} Assume that $V$ is  an irreducible representation of
$U_\eps$. Let  $\Phi_1,\Phi_2\in \Hom_{U_\eps}(V, V\o
U)$, and let $\chi_1, \chi_2\in \C[P]\o U[0]$ be the corresponding
generalized characters. Then

$$(\Phi_1, \Phi_2)=\frac{1}{\sqrt{\dim_\eps U}}(\chi_1, \chi_2)_1.\tag
4.7$$
\endproclaim
\demo{Proof} The proof repeats that of \cite{EK2} with minor changes
and is based on the identity $(\chi_1\o \bar \chi_2)_U=\chi_\Psi$,
where the intertwiner $\Psi: V\o V^*\to V\o V^*$ is given by
$$
\Psi=\vcenter{\epsfbox{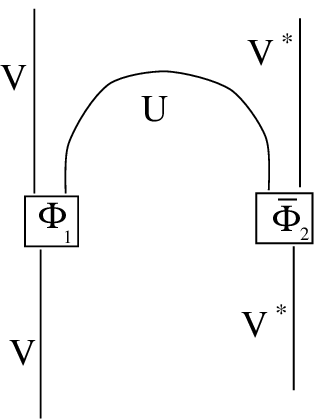}}
$$

This along with well-known identity $[\chi_\l \delta\bar
\delta]_0=\delta_{\l, 0}(-1)^{|R^+|}|W|$ proves the theorem.
\qed \enddemo

\head 5. Macdonald's theory.\endhead

In this section we consider an example of action of modular group in
modular tensor category obtained from quantum groups at roots of
unity. Namely,  we consider $\g=\sln$ and $U$ -- symmetric power of
fundamental representation. We will show that in this case the
$S$-matrix can be written in terms of Macdonald's polynomials of type
$A_{n-1}$ and deduce from this certain identities for values of these
polynomials at roots of unity.

Even though in this case one can write explicitly $\rho, h\v$,
we will use the general notations as far as
possible.

As in Section~3, we consider the reduced category $\Cal C(\sln,
\varkappa)$, based on representations of $U_\eps$, where $\eps=e^{\pi\i
/\varkappa}$. Let us fix a positive integer $k$ and
assume that $\varkappa$ has the following form:

$$\varkappa=K+kh\v, \quad K\in \Z_+.\tag 5.1$$

Define

$$C_K=\{\l\in P^+| \<\l, \theta\v\>\le K\}.\tag 5.2$$
Equivalently, we can rewrite this condition as follows: for any $\a\in
R^+$,

$$\<\l+k\rho, \a\v\><\varkappa-(k-1).\tag 5.3$$

Note that $C_K$ is non-empty and $\l\in C_K\iff \l^*\in
C_K$.

\example{Example} For $k=1$, this coincides with the domain $C$ we
defined in Section~3.
\endexample

Let $U=V_{(k-1)n\omega_1}$, where $\omega_1$ is the first fundamental
weight; in other words, $U$ is the deformation of the module
$S^{(k-1)n}\C^n$, where $\C^n$ is the fundamental representation of
$\sln$.  Note that due to  Lemma~3.2, $U$ is an irreducible
module over $U_\eps$. It will be extremely important for us that
$U[0]$ is one-dimensional; we fix some non-zero vector $u_0\in U[0]$,
which allows us to identify $U[0]\simeq \C: u_0\mapsto 1$.

\proclaim{Theorem 5.1} Let $\mu\in C$. Then

$$\dim \Hom_{\Cal C(\sln, \varkappa)} (V_\mu, V_\mu\o U)=
\cases
	1, \quad \mu=\lambda+(k-1)\rho \text{ for some }\lambda\in C_K\\
	0  \quad \text{ otherwise}  \endcases.$$

\endproclaim

Note that  $\l\in C_K$ implies $\l+(k-1)\rho\in  C$, so
$V_{\l+(k-1)\rho}$ is irreducible over $U_\eps$.

We will prove this theorem later.

{}From now on, let us for simplicity denote

$$\lk=\l+(k-1)\rho.$$

As before, let us denote $H=\bigoplus_{\mu\in C} V_\mu^*\o V_\mu$.
Theorem 5.1  allows us to choose a basis in $\Hom (H,
U)$. Indeed,   we have canonical isomorphism $\Hom(V_\mu^*\o
V_\mu, U)\simeq \Hom (V_\mu, V_\mu\o U)$. For $\l\in C_K$, let
$\Phi_\l:  V_\lk\to V_\lk\o U$ be an intertwiners such that
$\Phi(v_\lk)=v_\lk\o u_0+\dots$. It follows from Theorem~5.1 that such
an intertwiner exists and is unique and that

$$\Hom(H, U)\simeq \bigoplus_{\l\in C} \Hom (V_\l, V_\l\o U)
=\bigoplus_{\l\in C_K} \C \Phi_\l. $$

The main result of this section is that in this basis the action of
the matrix $S_U$ defined in Theorem~1.10 is given by the values of
Macdonald's polynomials at special points.

Recall that Macdonald's polynomials $P_\l^{q, q^k}$ (where $k$ is the
same positive integer that we used in the beginning of this section)
are elements of $\C(q)[P]^W$ which are defined by the following
conditions (see \cite{M1, M2}):

\roster\item
 $P_\l=e^\l+\text{\it lower
order terms}$
\item  $(P_\l, P_\mu)_k=0$ if $\l\ne\mu$, where

$$(f,g)_k=\frac{(-1)^k}{|W|} [f \bar
g\delta_k\overline{\delta_k}]_0. \tag 5.4$$
Here

$$\delta_k=\prod_{i=0}^{k-1}\prod_{\a\in R^+}
 (e^{\a/2}-q^{-2i}e^{-\a/2}),$$
and all other notations are as in  Section~4 with  complex
conjugation in $\C$ extended to $\C(q)$ by $\bar q=q^{-1}$.
\endroster

\remark{Remark} This definition, as well as the complex conjugation on
$\C(q)[P]$ differs from the definition in both  original Macdonald's
papers and  \cite{EK1--3}, which use the $\C(q)$-linear
inner product rather than hermitian. However, it is easy to check that
this definition is in fact equivalent to the original one, which
relies on the identity

$$P_\l^{q^{-1}, q^{-k}}=P_\l^{q, q^k}$$
(see \cite{M1}).

 We use the same notations as we did
in \cite{EK2}; thus, what we denote by $P_\l^{q, q^k}$ in the original
notations of Macdonald would be $P_\l(x; q^2, q^{2k})$.

\endremark

{}From now on, we will drop the superscript $q, q^{-k}$ and denote
Macdonald's polynomials simply  $P_\l$. The following properties of
these polynomials can be easily deduced from the definition:

$$\gathered
\overline{P_\l}=P_{\l^*},\\
\overline{P_\l(q^\mu)}=P_{\l^*} (q^\mu)=
P_\l(q^{-\mu})=P_\l(q^{\mu^*}).\endgathered\tag 5.5$$

Here $\overline{\phantom{T}}$ is the involution in $\C_q$ (in the
second line) and in $\C_q[P]$ (in the first line) which was defined in
Section~4.

Our arguments will be based on the relation between Macdonald's
polynomials of type $A$ and representations of $U_q \sln$. We recall
the main facts here, following the papers \cite{EK2, EK3}; note,
however, that the quantum group used in these papers differs from the
one used here  by  substitution  $q\leftrightarrow q^{-1}$.

For the
moment, we consider representations of $U_q \sln $ for generic $q$,
i.e. over the field $\C_q= \C(q^{1/2n})$. Let $k, U, u_0$ be the same as
above. Then we have the following results (see \cite{EK2, EK3}):

\roster
\item

$$\dim\Hom_{U_q\sln} (V_\mu, V_\mu\o U)=
\cases
	1, \quad \mu=\lambda+(k-1)\rho \text{ for some }\lambda\in P_+\\
	0  \quad \text{ otherwise}  \endcases.$$

We fix an intertwiner $\Phi_\l:V_\lk\to V_\lk\o U$ such that
$\Phi_\l v_\lk=v_\lk\o u_0 +\dots$

\item Let $\varphi_\l\in \C(q)[P]\o U[0]$ be the generalized character
of $\Phi_\l$ (see (3.14)).  Then

$$\gathered
P_\l=\frac{\varphi_\l}{\varphi_0},\\
\varphi_0=\prod_{\a\in
R^+}\prod_{i=1}^{k-1}(e^{\a/2}-q^{-2i}e^{-\a/2})\cdot u_0.\endgathered
\tag 5.6$$

\item

$$(\Phi_\l, \Phi_\l)=\frac{(-1)^{k-1}}{\sqrt{\dim_q U}}
(u_0, u_0)(P_\l, P_\l)_k,$$
where
$$(P_\l, P_\l)_k=
\prod_{\a\in R^+}\prod_{i=1}^{k-1}
\frac{[(\a, \l+k\rho)+i]}
     {[(\a, \l+k\rho)-i]}\tag 5.7$$

This identity is a reformulation of famous Macdonald's inner product
identities in our case, i.e. for hermitian rather than bilinear inner
product.
\endroster

Now we can come back to Theorem~5.1.

\demo{Proof of Theorem 5.1} It follows from Lemma~3.7 that $\dim
\Hom_{\Cal C} (V_\mu,V_\mu\o U)\le 1$, and it can be non-zero only
if $\mu=\l+(k-1)\rho, \l\in P^+$. It follows from Lemma~4.3 that this
dimension is equal to one iff
$(\Phi_\l, \Phi_\l)\ne 0$, where $\Phi_\l$ is the corresponding
intertwiner for $U_\eps$. On the other hand, formula (5.7) shows that
$(\Phi_\l, \Phi_\l)\ne 0\iff \l\in C_K$. \qed\enddemo

Note that the proof used highly non-trivial result -- explicit formula
for the norm $(\Phi_\l, \Phi_\l)$(Macdonald's inner product formula).

Now, let us come back to the case of roots of unity. Let $\l\in C_K$,
and let $\Phi_\l\in \Hom_{\Cal C(\g, \varkappa)}(V^*_\lk\o V_\lk,  U)$
be as before.

\proclaim{Theorem 5.2} Let $\l\in C_K$. Then $P_\l$ is well defined at
$q=\eps$ \rom{(}i.e., its coefficients, which are rational functions of $q$,
are well-defined at $q=\eps$\rom{)}. \endproclaim

\demo{Proof} This follows from the fact that Macdonald polynomials can
be written in terms of generalized characters (see (5.6) above) and
Lemma~3.7.\qed\enddemo

Now comes the crucial step.

\proclaim{Lemma 5.3} We have the following identity in the category
$\Cal C(\sln, \varkappa)$:

$$\vcenter{\epsfbox{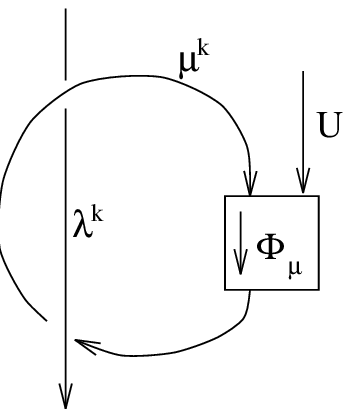}}
=\varphi_\mu^\eps(\eps^{-2(\l+k\rho)})  \Phi_\l, $$
where $\varphi_\l^\eps$ is the element of $\C[P]$ which is obtained by
substituting  $q=\eps$ in the expression for $\varphi_\l$ and
identifying $U[0]\simeq \C: u_0\mapsto 1$.
\endproclaim

\demo{Proof} For the case when $q$ is indeterminate, it was proved in
\cite{EK3}; it is easy to see that in fact all the arguments can be
carried out in the case $q=\eps$ as well.
\qed\enddemo

This immediately implies the following theorem:

\proclaim{Theorem 5.4} Let $\g=\sln$, and let $U, \Phi_\l, \l \in C_K$
be as above.
Then the action of the modular group in $\Hom(H, U)$ in this basis is
given by the matrices $S_U=(S_{\l\mu}), T_U=(T_{\l\mu})$, where

$$\gathered
T_{\l\mu}=\delta_{\l\mu} \eps^{(\l+k\rho, \l+k\rho)
	-\frac{\varkappa}{n}(\rho,\rho)},\\
S_{\l\mu}=d_\l  P^\eps_\mu(\eps^{-2(\l+k\rho)}), \endgathered\tag 5.8$$
where $P^\eps_\l$ is
Macdonald's polynomial $P_\l^{q, q^k}$ calculated  at $q=\eps$, and

$$d_\l = \frac{\i^{n(n-1)/2}}{\sqrt n \varkappa^{(n-1)/2}}
\prod_{\a\in R^+}\prod_{i=0}^{k-1} \bigl(\eps^{-(\a,
\l+k\rho)}-\eps^{-2i+(\a, \l+k\rho)}\bigr).\tag 5.9$$
\endproclaim

\demo{Proof} Formula for $T$ is obvious from formula (3.19) for $\zeta$
 and $(\lk, \lk+2\rho)= (\l+k\rho, \l+k\rho)-(\rho, \rho)$. It follows
from the definition of $S$ and Lemma~5.3 that $S_{\l\mu}$ is given by
formula (5.8) with

$$d_\l=\frac{\dim_\eps V_{\l^k}}{D}\varphi_0
(\eps^{-2(\l+k\rho)})$$
(as before, we consider $\varphi_0$ as scalar-valued).
Substituting in this expression Weyl formula for $\dim_\eps V_{\l^k}$,
expression (3.19) for $D$ and formula (5.6) for $\varphi_0$, we see
that
$$d_\l= \frac{\i^{|R^+|}}{\sqrt{|P/\varkappa Q|}}
\prod_{\a\in R^+}\prod_{i=0}^{k-1}
\bigl(\eps^{-(\a, \l+k\rho)}-\eps^{-2i+(\a, \l+k\rho)}\bigr).$$
Since for $\sln$ we have $|P/Q|=n, |R^+|=n(n-1)/2$, and the rank  is
 $n-1$, we get formula (5.9).\qed\enddemo

Similar formulas for the action of $SL_2(\Z)$ in terms of the values of
Macdonald's polynomials were obtained by Cherednik (\cite{Ch}) in the
study of difference Fourier transform.

\example{Example} Consider the case $\g=\sltwo$. Then every
irreducible finite-dimensional representation has the form
$V_{(k-1)n\omega_1}$  for
some choice of $k$, and thus, in this case Theorem~5.4 gives all
matrix coefficients of the action of the modular group in $H$, which
in this case are written in terms of $q$-ultraspherical polynomials
(=Macdonald's polynomials for $\sltwo$). In particular, this shows
that $S$-matrix can be written in terms  of
basic hypergeometric functions with parameter $q$ taken to be root of
unity (see \cite{AI} for expressions of $q$-ultraspherical polynomials in
terms of basic hypergeometric functions).
\endexample

Explicit calculation, using symmetry properties (5.5) of Macdonald's
polynomials, gives the following symmetries of $S$-matrix:

$$\gathered
S_{\l\mu}=S_{\l^*\mu^*},\\
\overline{S_{\l\mu}}=(-1)^{(k-1)n(n-1)/2}\eps^{n(n-1)k(k-1)/2} S_{\l^*\mu}.
\endgathered\tag 5.10$$

Also, it is easy to calculate the action of the matrix $C$ in the
basis $\Phi_\l$. As before, let us assume that we have chosen
identifications $V_{\l^*}\simeq V_\l^*$ as in (1.1). Then
we have the following theorem:

\proclaim{Theorem 5.5}

$$\Phi_\l\ C^{-1} =(-1)^{(k-1)n(n-1)/2}
		\eps ^{n(n-1)k(k-1)/2}\Phi_{\l^*}.\tag 5.11$$

\endproclaim
\demo{Proof} Let $v_\lk$ be highest-weight vector in $V_\lk$,
$v_\lk^*$ -- lowest weight vector in $V_{\lk}^*$ such that $\<v_\lk^*,
v_\lk\>=1$. Also, let $w_\lk$ be lowest weight vector in $V_\lk$,
$w_\lk^*$ -- highest weight vector in $V_{\lk}^*$ such that $\<w_\lk^*,
w_\lk\>=1$. Then by definition $\Phi_\l (v_\lk^*\o v_\lk)=u_0$, and it
follows from (5.6) and symmetry of Macdonald's polynomials that
$\Phi_\l (w_\lk^*\o w_\lk)=u_0 (-1)^{(k-1)n(n-1)/2}\eps ^{-n(n-1)k(k-1)/2}$.

It follows from formula (3.6) for universal $R$-matrix that

$$\aligned
\Phi_\l C^{-1}:&V_\l \o V_\l^*\to U\\
& v_\lk\o v_\lk^*\mapsto \theta_\lk \eps ^{-(\lk, \lk)}
\Phi_\l (v_\lk^*\o v_\lk)= \eps^{(\lk, 2\rho)} u_0.\endaligned$$

On the other hand, if we identify $V_\l\simeq (V_{\l^*})^*,
V_\l^*\simeq V_{\l^*}$ and denote by $\<\, , \,\>$ canonical pairing
$(V_{\l^*})^*\o V_{\l^*}\to \C$ then

$$\<v_\lk, v_\lk^*\>=\eps^{(\lk, 2\rho)}.$$
Therefore, similarly to what we discussed before,

$$\Phi_{\l^*} (v_\lk\o v_\lk^*)= (-1)^{(k-1)n(n-1)/2}\eps^{-n(n-1)k(k-1)/2}
\eps ^{(\lk, 2\rho)} u_0.$$

Comparing these expressions, we get the statement of the theorem.
\qed\enddemo

\remark{Remark} Note that since $\theta_U=\eps^{n(n-1)k(k-1)}$, which
is verified by direct computation, we again see that
$C^2=\theta_U^{-1}$. \endremark

Now we can rewrite results about the action of modular group which
were proved in purely abstract setting in Sections~1 and 2 to this case,
which results in identities for Macdonald's polynomials:

\proclaim{Theorem 5.6} For $\l, \mu\in C_K$,

$$S_{\l\mu}(P_\l, P_\l)_k=S_{\mu\l} (P_\mu, P_\mu)_k. \tag 5.12$$

\endproclaim

\demo{Proof} This is nothing but the condition of unitarity of matrix $S$
with respect to the inner product on intertwiners
(Theorem~2.5). Indeed, the unitarity condition can be rewritten as
follows:
$(\Phi_\mu S, \Phi_\l)=(\Phi_\mu, \Phi_\l C^{-1} S),$
which is equivalent to

$$S_{\l\mu}(\Phi_\l,\Phi_\l)=
(-1)^{(k-1)n(n-1)/2}\eps^{-n(n-1)k(k-1)/2}\overline{S_{\mu\l^*}}
(\Phi_\mu, \Phi_\mu).$$

Using  symmetry properties (5.10), we get the statement of the theorem.
\qed\enddemo

Using expressions for $S$-matrix given in Theorem~5.4, we can rewrite
this result as follows:

$$P_\l^\eps(\eps^{-2(\mu+k\rho)}) (P_\mu, P_\mu)_k d_\mu=
P_\mu^\eps(\eps^{-2(\l+k\rho)}) (P_\l, P_\l)_k d_\l,\tag 5.13$$
which is precisely the symmetry identity for Macdonald's
polynomials of type $A$ (see \cite{EK3}).  One can check that in fact
all our arguments work for general $q$, i.e. one can avoid using the
fact that the category is modular; essentially, this is the same proof
that was given in \cite{EK3}, only now it has a clear interpretation.

Other identities, which only apply to modular categories and thus
do not generalize to the case of indeterminate $q$  can be obtained
from the relations in modular group. This gives the following purely
combinatorial theorem:

\proclaim{Theorem 5.7} Let $S=(S_{\l\mu}), T=(T_{\l\mu}), \l,\mu\in
C_K$ be the matrices  given by \rom{(5.8), (5.9)}. Then
$$\gathered
S^2=(-1)^{(k-1)n(n-1)/2}\eps^{-n(n-1)k(k-1)/2}\delta_{\l\mu^*},\\
(ST)^3=S^2.\endgathered
\tag 5.11$$
\endproclaim

These are  certain identities   for Macdonald's polynomials at roots of
unity, which were not known before and which would be very difficult to
prove by combinatorial methods.
 Again, similar (and even more general)
identities have been recently obtained by
Cherednik (\cite{Ch}) in the study of difference Fourier transform
related with double affine Hecke algebras.

\head 6. Characters and Grothendieck ring\endhead

In this section we again return to consideration of the category $\Cal
C(\g, \varkappa)$ for arbitrary $\g$ and describe its Grothendieck ring.
We also give an elementary proof of the fact that the matrix $s_{\l\mu}$
defined by (3.18) is non-degenerate, and calculate its square.
Results of this section are not new, but I
was unable to locate them in the literature\footnote{For example, it
is mentioned as unknown for exceptional Lie algebras in \cite{Tu,
XI.6.4}}, so for the sake of
completeness they are included here.

Recall (see Section~3) that  we have fixed $\varkappa\in \Z_+$ such
that $\varkappa\ge  h\v$, and we have defined the open and closed alcoves

$$\gathered
C=\{\l\in P^+|\<\l+\rho, \theta\v\><\varkappa\},\\
\overline{C}=\{\l\in P|\<\l+\rho, \a_i\v\>\ge 0,
\<\l+\rho, \theta\v\>\le\varkappa\}. \endgathered
$$

We have also defined affine Weyl group
$\What=W\ltimes \varkappa Q\v$ and its shifted action on $\h^*$ by
$w.\l=w(\l+\rho)-\rho$.
Then, as is well-known, we have the following statements:
\roster
\item $\overline{C}$ is the fundamental domain for the shifted action of
$\What$.

\item Every $\l\in C$ is regular with respect to the shifted action of
$\What$: $w.\l=\l$ iff $w=1$.

\item  For every $f\in \C[P]^W, \mu \in P, w\in \What$ we have
$$f(\eps ^{2w(\mu)})=f(\eps^{2\mu}).$$

\endroster

It fact, the last statement can be reversed: if $\l, \mu\in P$ are
such that  $f(\eps^{2\l})=f(\eps^{2\mu})$ for all $f\in \C[P]^W$ then
$\l=w(\mu)$ for some $w\in \What$; however, we won't use this result.

For every $\l\in P$ define $\chi_\l\in \C[P]^W$ by

$$\chi_\l=\frac{\sum_W (-1)^{l(w)}e^{w(\l+\rho)}}{\delta},\tag 6.1$$
where $\delta$ is Weyl denominator (3.9).
For $\l\in P^+$, $\chi_\l$ is the character of the module $V_\l$.

Now let $\eps=e^{\pi\i/m\varkappa}$. Recall that we denote

$$\dim_\eps V_\l= \Tr_{V_\l} (\eps^{2\rho})
= \chi_\l(\eps^{2\rho})=\chi_\l(\eps^{-2\rho}).$$

 It is easy to see  that $\dim_\eps V_\l=
0$ for $\l\in \bar C\setminus C$, and $\dim_\eps V_\l\ne 0$ for $\l\in C$.

For $\l, \mu\in P$, define the numbers  $s_{\l\mu}\in \C$ by

$$s_{\l\mu}=
\frac{\sum_{w\in W} (-1)^{l(w)}\eps^{-2(w(\l+\rho), \mu+\rho)'}}
	{\delta(\eps^{-2\rho})}.\tag 6.2
$$

If $\l, \mu\in C$ this can also be rewritten as follows:
$$s_{\l\mu}= \chi_\l (\eps^{-2(\mu+\rho)})\dim_\eps V_\mu.$$

\proclaim{Lemma 6.1}

$$\gathered
s_{\l\mu}=s_{\mu\l},\\
s_{\l\mu}=(-1)^{l(w)}s_{\l\ w.\mu}\quad\text{ for any }w\in\What,\\
s_{\l\mu}=0\quad \text{ if }\l\in \bar C\setminus C.\endgathered\tag 6.3
 $$

\endproclaim

\proclaim{Theorem 6.2} Let $s=(s_{\l\mu})_{\l, \mu\in C}$. Then

$$s^2=D^2c,$$
where $c_{\l\nu}=\delta_{\l\nu^*}$ and
$$D^2=|P/\varkappa  Q\v|\prod_{\a\in R^+} \left(2\sin\frac{(\a,
\rho)}{\varkappa}\pi \right)^{-2}
=|P/\varkappa  Q\v|(-1)^{|R^+|}\delta^{-2}(\eps^{-2\rho}).
\tag 6.4$$
\endproclaim

\demo{Proof} Let $\l, \mu\in C$. Then

$$\aligned
\sum_{\mu\in C}s_{\l\mu}s_{\mu\nu}=&\sum_{\mu\in \bar
C}s_{\l\mu}s_{\mu\nu}= \sum_{\mu\in P/\What}s_{\l\mu}s_{\mu\nu}=
\frac{1}{|W|} \sum_{\mu\in P/\varkappa Q\v}s_{\l\mu}s_{\mu\nu}\\
=&\frac{1}{|W| \delta^2(\eps^{-2\rho})} \sum_{\mu\in P/\varkappa Q\v}
\sum_{w, w'\in W} (-1)^{l(ww')}\eps^{-2(\mu+\rho,
w(\l+\rho)+w'(\nu+\rho))'}.\endaligned$$

For any $a\in P$ we have

$$\sum_{\mu\in P/\varkappa Q\v}\eps^{2(\mu+\rho, a)'}=
\cases 0, \quad a\notin \varkappa Q\v\\
       |P/\varkappa Q\v|, \quad a\in \varkappa Q\v.\endcases
$$

Since  $\l, \nu\in C$, it follows from the fact that $C$ is fundamental
domain for the action of $\What$  that
$w(\l+\rho)+w'(\nu+\rho)\in \varkappa Q\v$
is only possible if $\l=\nu^*, ww'=w_0$ --
the longest element in $W$. Thus,

$$\sum_{\mu\in C}s_{\l\mu}s_{\mu\nu}=\delta_{\l\nu^*} |P/\varkappa Q\v|
(-1)^{|R^+|}\delta^{-2}(\eps^{-2\rho}).$$
\qed \enddemo

\proclaim{Corollary 6.3} The matrix $s$ defined by \rom{(6.2)}
 is non-degenerate.\endproclaim

In a similar way, one can prove the identity $D\zeta^3 t^{-1}st^{-1}=
sts$, where $D$ is as above and $\zeta$ is given by (3.19) -- this
requires calculation of $\sum_\mu \eps^{2(\mu, \mu+a)}$.

Now denote by $K$ the Grothendieck ring of the category $\Cal C(\g,
\varkappa)$, and by $K_\C=K\o_\Z\C$ its complexification.

Let $F(C)$ be the ring of all complex-valued functions on $C$. For every
$V\in Rep\ U_\eps$ denote by $f_V\in F(C)$ the function given by

$$f_V(\mu)=\ch V( \eps^{-2(\mu+\rho)}).$$

\proclaim{Lemma 6.4} If $V$ is negligible then $f_V=0$ on $C$.\endproclaim

\demo{Proof} It follows from the following identity

$$f_V(\mu)=\frac{1}{\dim_\eps V_\mu}\quad\vcenter{\epsfbox{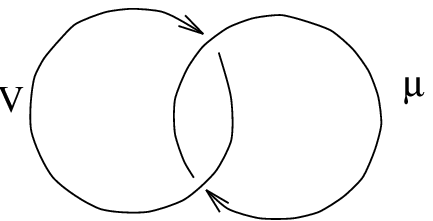}} $$

\enddemo

\proclaim{Corollary} The map $V\mapsto f_V$ is a well-defined ring
homomorphism
$K\to F(C)$.\endproclaim

\proclaim{Theorem 6.5} The map $V\mapsto f_V$ is an isomorphism
$K_\C\simeq F(C)$. \endproclaim

\demo{Proof} It suffices to prove that $\text{det }
(\chi_\l(\eps^{-2(\l+\rho)}))_{\l, \mu\in C}\ne 0$, which follows from
non-degeneracy of matrix $s$ (Corollary 6.3). \qed\enddemo

\proclaim{Corollary 6.6}

$$K_\C\simeq \C[P]^W/\Cal I,$$
where the ideal $\Cal I$ is spanned as a vector space
 by the elements of the form
$\chi_\l+ \chi_{s_{\Gamma}\l}, \l\in P$, where $s_{\Gamma}$ is the
reflection with respect to the affine wall $\Gamma$ \rom{(}see
\rom{(3.15))}.
\endproclaim

\remark{Remark} It can be shown (see  \cite{F}) that for $\varkappa$
large enough there is a stronger result: $K_\C\simeq \C[P]^W/\Cal
I$, and the ideal $\Cal I$ is generated as an ideal by $\chi_\l, \l\in
\Gamma\cap P^+$. In particular, this is always so for
$\g=\frak{sl}_n$. \endremark

\head 7. More on hermitian structure in $\Cal C (\g,
\varkappa)$.\endhead

In this section we give another description of the hermitian structure
 in $\Cal C (\g,\varkappa)$. It will be used in the future papers to
establish relation with affine Lie algebras. Also, it allows us to
define the inner product on the spaces of intertwiners uniquely up to
a constant from $\R_+$ (not from $\C^\times$, as we did in Section~4);
thus, it makes sense to discuss whether this inner product is positive
definite.

Recall that the key ingredient of the definition of hermitian structure
was definition of an involution $\omega$ on $\Ug$, which allowed us to
define for every module $V$ a module $V^\omega$, and isomorphisms
$V^\omega\simeq V^*$.

Here is another description of the same involution. As before, we begin with
consideration of generic $q$, i.e. of the quantum group over the field
$\C_q$.

Let the compact involution $\omega_c$ be the  antilinear
algebra automorphism $\Ug\to \Ug$ defined by

$$\aligned
	\omega_c:         & e_i\mapsto -q^{d_i}f_i,\\
		       & f_i\mapsto -q^{-d_i}e_i,\\
		       &q^h\mapsto q^{h},\\
	               & q\mapsto q^{-1}.
\endaligned\tag 7.1$$

One easily checks that $\omega_c$ is also coalgebra automorphism.
Similarly to the constructions of Section~4, for every $\Ug$-module $V$ and
homomorphism $\Phi$ we can define $V^{\omega_c}$ and
$\Phi^{\omega_c}$. Again, for an irreducible highest-weight module
$V_\l$ we have a (not canonical)  isomorphism  $V_\l^{\omega_c}\simeq
V_{\l}^*$; due to complete reducibility, the same is true for
arbitrary module $V$.

However, this involution can not be used to define a hermitian
structure on the category of representations because there is no
canonical isomorphism between $(V\o W)^{\omega_c}$ and $W^{\omega_c}\o
V^{\omega_c}$; instead, we have isomorphism $(V\o W)^{\omega_c}\simeq
V^{\omega_c}\o W^{\omega_c}$.

To get a hermitian structure, we need one more ingredient, namely, the
longest element of the quantum Weyl group, which was studied in
\cite{LS, L7}. We reformulate the results of these papers in the
following theorem:

\proclaim{Theorem 7.1}\rom{(Levendorskii--Soibelman, Lusztig)}
There exists an element $\Omega$ in a certain completion of $\Ug$
satisfying the following properties:

\roster\item
$\Omega$ acts in every finite-dimensional $\Ug$-module, and
$\Omega V[\l] \subset V[w_0(\l)]$.

\item
$$\gathered
\Omega f_i \Omega^{-1}= - q^{-d_i} e_{i\v}\\
\Omega e_i \Omega^{-1}= - q^{d_i} f_{i\v}\\
\Omega q^h \Omega^{-1}= q^{w_0(h)}.\endgathered \tag 7.2
$$
\item $$\Delta(\Omega)=\Cal R^{-1}(\Omega\o \Omega)$$

\item $\Omega^2 = Z\theta^{-1}$, where $\theta$ is the universal twist
shown on Fig.~1A \rom{(}recall that $\theta$ is a central element such
that  $\theta|_{V_\l}=q^{(\l,
\l+2\rho)'}\Id$\rom{)}, and $Z$ is a central element satisfying $\Delta
Z=Z\o Z$.
\endroster\endproclaim

\remark{Remark 7.2} It is convenient to think of $\Omega$ as
represented by the following ribbon graph:

$$
\vcenter{\epsfbox{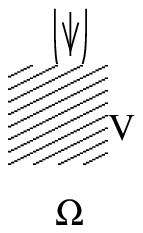}}
$$

Such a graph is not allowed in the original formalism developed by
Reshetikhin and Turaev (and indeed, $\Omega$ is not an intertwining
operator); however, this becomes possible after suitable extension of
the formalism. \endremark

Comparing (7.2) with definition of  $\omega$ in Lemma~4.1 we see that

$$\omega (x)= \omega_c\left(\Omega x \Omega^{-1}\right).\tag 7.3$$

Now we can come back to defining the isomorphisms $V^\omega\simeq
V^*$.  Such an isomorphism is equivalent to defining
a pairing $V^\omega\o V\to \C_q$, or a
non-degenerate hermitian form on $V$ such that

$$H(xv, v')=H(v, S\omega(x) v').\tag 7.4$$

Similarly, isomorphism $V^{\omega_c}\simeq V^*$ is equivalent to defining a
non-degenerate hermitian form $(\, , \,)_c$ on $V$ satisfying

$$(xv, v')_c=(v, S\omega_c(x) v')_c.\tag 7.5$$

In fact, given any of these forms we can define another:

\proclaim{Lemma 7.3} If a hermitian form $(\, , \,)_c$ on $V$
satisfies condition \rom{(7.5)} then the form

$$H(v, v')=(\Omega v,  v')_c\tag 7.6$$
satisfies condition \rom{(7.4)} and vice versa. \endproclaim

\demo{Proof} Obvious from (7.3).\qed\enddemo

So far, we have considered the case of indeterminate $q$. Now, let us
specify $q=\eps=e^{\pi\i/m\varkappa}$.

\proclaim{Theorem 7.4} Let $V_\l, \l\in C$ be an irreducible  highest weight
module over $U_\eps$, and $v_\l$ -- highest weight vector. Let $(\, ,
\,)_c$ be the hermitian form on $V_\l$ satisfying \rom{(7.5)} and
normalized  by condition $(v_\l, v_\l)_c=1$. Then this form is positive
definite. \endproclaim

\demo{Proof} For $q=1$ this is well-known. Now, let
$\eps^t=e^{t\pi\i/m\varkappa}, t\in [0,1]$. For every $t$ we can
define a module $V_\l$ over $U_{\eps^t}$ and a form $(\, , \,)_c^t$ as
in the theorem. We can identify all $V_\l$ (as vector spaces over
$\C$) and thus get a family of forms on the same space.   It is easy
to check, using results of Section~3(see (3.21), (3.22)),
 that for every $t\in [0,1]$ the
module $V_\l$ is irreducible over $U_{\eps^t}$, and thus the form $(\,
, \,)_c^t$ is non-degenerate. Since for $t=0$ this form is
positive-definite, the same holds for all $t$, in particular, for
$t=1$. \qed\enddemo

\remark{Remark 7.5} This does not hold in more general situation of
rational $\varkappa$ (cf. Remark~3.10). \endremark

Therefore, we arrive to the following theorem, which is the main
result of this section:

\proclaim{Theorem 7.6} In the notations of Theorem~\rom{7.4}, there
exists a unique up to a positive real constant hermitian form $H$ on
$V_\l$ which satisfies the invariance condition \rom{(7.4)} and such
that the form
$(v, v')_c=H(\Omega^{-1}v, v')$ is positive definite. This form $H$
can be defined by the condition $H(\Omega^{-1}v_\l, v_\l)\in
\R_+$.\endproclaim

This defines the form $H$ (and thus, isomorphism $V^\omega\simeq V^*$)
on irreducible modules. We can extend it to tensor products by the
rule $H(v\o w, v'\o w')=H(v, v')H(w, w')$ (this satisfies the
invariance condition due to the fact that $\omega$ is coalgebra
antiautomorphism). Therefore, we can define the inner product on every
space of intertwiners $\Hom_{U_\eps} (V_\l, V_\mu\o V_\nu)$ uniquely
up to a real positive factor.

\proclaim{Conjecture 7.7} So defined inner product on
$\Hom_{\Cal C(\g, \varkappa)} (V_\l, V_\mu\o V_\nu)$ is positive
definite.
\endproclaim

In the simplest case $\g=\sltwo$ this can be checked directly. In
general case, the answer is not known.


\Refs
\widestnumber\key{AAA}

\ref\key A\by Andersen, H. H. \paper On tensor products of quantized
tilting modules \jour Com. Math. Phys.\vol 149\yr 1992\pages
149--159\endref



\ref\key AI \by Askey, R. and Ismail, M.E.H.
\paper A generalization of ultraspherical polynomials
\inbook Studies in Pure Mathematics
\ed P. Erd\"os\publ Birkh\"auser \yr 1982\pages 55--78\endref


\ref\key AP\by Andersen, H. H.  and Paradowski, J.
\paper Fusion categories arising from semisimple Lie algebras
\jour Com. Math. Phys\vol 169\yr 1995\pages 563--588\endref

\ref\key Ch\by Cherednik, I.\paper Macdonald's evaluation conjectures
and difference Fourier  transform
\jour pre\-print, May  1995, q-alg/9412016\endref



\ref\key Dr1\by Drinfeld, V.G. \paper Quantum groups\inbook Proc. Int.
Congr. Math., Berkeley, 1986\pages 798--820\endref

\ref\key Dr2 \bysame\paper On almost cocommutative Hopf
algebras \jour Leningrad Math.J. \vol 1\issue 2\yr 1990\pages
321--342\endref

\ref\key EK1\by Etingof, P.I. and Kirillov, A.A., Jr\paper A unified
representation-theoretic approach to special functions
\jour  Functional Anal. and its Applic.\vol 28\issue 1
 \yr 1994\pages 91--94\endref

\ref\key EK2\bysame\paper Macdonald's
polynomials and representations of quantum groups \jour
Math. Res. Let.\vol 1\yr 1994\pages 279--296\endref


\ref \key EK3 \bysame\paper
Representation-theoretic proof of inner product and symmetry
identities for Macdonald's polynomials\jour hep-th/9410169, to appear in
 Comp. Math. \yr 1995\endref

\ref\key F\by Finkelberg, M.
\paper Fusion categories\jour Ph.D thesis, Harvard Univ.\yr
1993\finalinfo (to appear in GAFA)\endref


\ref\key GK\by Gelfand, S. and Kazhdan, D.
\paper Examples of tensor categories \jour Invent. Math. \vol 109\yr
1992\pages 595--617\endref



\ref\key Kac\by Kac, V.G. \book Infinite-dimensional Lie algebras\publ
Cambridge Univ. Press\bookinfo 3rd ed.\yr 1990\endref

\ref\key Kas\by Kassel, C.\book Quantum groups\publ Springer\publaddr
New York\yr 1995\endref


\ref\key KL1\by Kazhdan, D. and Lusztig, G.\paper Tensor structures
arising from affine Lie algebras. {\rm I} \jour J. of
AMS\vol 6\yr 1993\pages 905--947\endref

\ref\key KL2\bysame\paper Tensor structures
arising from affine Lie algebras. {\rm II} \jour J. of
AMS\vol 6\yr 1993\pages 949--1011\endref

\ref\key KL3\bysame\paper Tensor structures
arising from affine Lie algebras. {\rm III} \jour J. of
AMS\vol 7\yr 1994\pages 335--381\endref

\ref\key KL4\bysame\paper Tensor structures
arising from affine Lie algebras. {\rm IV} \jour J. of
AMS\vol 7\yr 1994\pages 383--453\endref


\ref\key Ke \by Kerler, T.\paper Mapping class group actions on
quantum doubles\jour Comm. Math. Phys.\vol 168 \yr
1995\pages 353--388\endref



\ref\key L1\by Lusztig, G.\paper Modular representations and quantum
groups\inbook Contemp. Math.\vol 82\yr 1989\publ AMS\publaddr
Providence\pages 52--77\endref

\ref \key L2\bysame \paper Finite-dimensional Hopf algebras arising
from quantized universal enveloping algebras\jour J. of AMS\vol 3\yr
1990\pages 257--296\endref

\ref \key L3\bysame \paper Quantum groups at roots of 1\jour
Geom. Dedicat. \vol 35\yr 1990\pages 89--113\endref

\ref \key L4\bysame \paper On quantum groups \jour J. of Algebra\vol 131\yr
1990\pages 466--475\endref



\ref\key L5\bysame\book Introduction to quantum groups\publ
Birkh\"auser \publaddr Boston \yr 1993\endref

\ref \key L6\bysame \paper Monodromic systems on affine flag manifolds
\jour Proc. R. Soc. Lond. A \vol 445\yr 1994\pages 231--246\endref

\ref \key L7\bysame \paper Canonical bases arising from quantized
 enveloping algebras. {\rm II}\jour Progr. Theor. Phys. Suppl. \vol
102\yr 1990 \pages 175--201\endref

\ref\key LS\by Levendorskii, S.Z. and Soibelman, Ya. S. \paper Some
applications of quantum Weyl groups\jour Jour. Geom. Phys.\vol 7\issue
2\yr 1990\pages 241--254\endref

\ref\key Lyu\by Lyubashenko, V.\paper Modular transformations for
tensor categories \jour Preprint \yr 1993\endref



\ref\key M1\by Macdonald, I.G. \paper A new class of symmetric
functions\jour Publ. I.R.M.A. Strasbourg, 372/S-20, Actes 20
S\'eminaire Lotharingien\pages 131-171\yr 1988\endref


\ref\key M2\bysame\paper Orthogonal polynomials associated
with root systems\jour preprint\yr 1988\endref


\ref\key Mac \by MacLane, S. \book Categories for working
mathematician \bookinfo Graduate Texts in Mathematics\vol 5 \publ
Springer--Verlag \publaddr New York\yr 1971\endref

\ref\key MS1\by Moore, G., Seiberg, N. \paper Classical and quantum
conformal field theory\jour Com. Math. Phys.\vol 123 \pages
177--254\yr 1989\endref



\ref\key MS2 \bysame \paper Lectures on RCFT
\inbook Superstrings '89 (Proc. of the 1989 Trieste Spring   School)
\eds M.~Green et al
\publ World Sci.
\publaddr River Edge, NJ
\yr 1990
\pages 1--129
\endref

\ref\key MS3\bysame \paper Polynomial equations for rational conformal
field theories\jour Phys. Let. \vol B212\yr 1988\pages 451--460\endref

\ref\key MS4\bysame \paper Naturality in conformal field theory\jour
Nucl. Phys.\vol B313\yr 1988 \pages 16--40\endref


\ref\key RT1 \by Reshetikhin, N. and Turaev, V. \paper
Ribbon graphs and their invariants derived from quantum groups
\jour Comm. Math. Phys.\vol 127 \pages 1--26\yr 1990\endref

\ref\key RT2 \by Reshetikhin, N. and Turaev, V. \paper Invariants of
3-manifolds via link polynomials and quantum groups\jour Inv. Math.
\vol 103\pages 547-597\yr 1991\endref


\ref \key Tu \by Turaev, V.\book Quantum invariants of knots and
3-manifolds \publ W. de Gruyter\publaddr Berlin\yr 1994\endref

\ref\key Vaf\by Vafa, C.\paper Towards classification of conformal theories
\jour Phys. Lett. B\vol 206 \yr 1988\pages 421--426\endref

\ref\key We\by Wenzl, H.\paper Quantum groups and subfactors of type
B, C, and D \jour Comm. Math. Phys.\vol 133\yr 1990\pages
383--433\endref

\endRefs
\enddocument
\end